\documentclass[12pt,oneside,reqno]{amsart}
\usepackage{txfonts}
\usepackage{bbm}
\usepackage{amsmath}
\usepackage{graphicx}
\usepackage{mathrsfs}
\usepackage{stmaryrd}
\usepackage{amsfonts}
\usepackage{enumerate,amsmath,amssymb,amsthm}

\numberwithin{equation}{section}

\newcommand{\be}{\begin{eqnarray}}
\newcommand{\ee}{\end{eqnarray}}
\newcommand{\ce}{\begin{eqnarray*}}
\newcommand{\de}{\end{eqnarray*}}
\newtheorem{theorem}{Theorem}[section]
\newtheorem{lemma}[theorem]{Lemma}
\newtheorem{remark}[theorem]{Remark}
\newtheorem{definition}[theorem]{Definition}
\newtheorem{proposition}[theorem]{Proposition}
\newtheorem{Examples}[theorem]{Example}
\newtheorem{corollary}[theorem]{Corollary}

\def\eps{\varepsilon}

\def\a{\alpha}

\def\b{\beta}
\def\p{\partial}

\def\[{{\Big[}}
\def\]{{\Big]}}
\def\<{{\langle}}
\def\>{{\rangle}}
\def\({{\Big(}}
\def\){{\Big)}}

\def\bx{{\mathbf{x}}}

\def\dif{{\mathord{{\rm d}}}}

\def\no{\nonumber}
\def\={&\!\!=\!\!&}
\def\bt{\begin{theorem}}
\def\et{\end{theorem}}
\def\bl{\begin{lemma}}
\def\el{\end{lemma}}
\def\br{\begin{remark}}
\def\er{\end{remark}}

\def\bd{\begin{definition}}
\def\ed{\end{definition}}
\def\bp{\begin{proposition}}
\def\ep{\end{proposition}}
\def\bc{\begin{corollary}}
\def\ec{\end{corollary}}
\def\bx{\begin{Examples}}
\def\ex{\end{Examples}}
\def\cA{{\mathcal A}}
\def\cB{{\mathcal B}}

\def\cD{{\mathcal D}}

\def\cJ{{\mathcal J}}

\def\cM{{\mathcal M}}

\def\mE{{\mathbb E}}

\def\mH{{\mathbb H}}

\def\mN{{\mathbb N}}

\def\mP{{\mathbb P}}

\def\mR{{\mathbb R}}
\def\mS{{\mathbb S}}

\def\sA{{\mathscr A}}

\def\sC{{\mathscr C}}

\def\sF{{\mathscr F}}

\def\sL{{\mathscr L}}
\def\sM{{\mathscr M}}

\def\geq{\geqslant}
\def\leq{\leqslant}

\def\div{\mathord{{\rm div}}}

\allowdisplaybreaks

\begin{document}

\title{Well-posedness and large deviation for degenerate SDEs
with Sobolev coefficients}

\date{}
\author{Xicheng Zhang}

\dedicatory{
School of Mathematics and Statistics,
Wuhan University, Wuhan, Hubei 430072, P.R.China\\
email: XichengZhang@gmail.com
 }


\begin{abstract}
In this article we prove the existence and uniqueness for degenerate stochastic differential equations
with Sobolev (possibly singular) drift and diffusion coefficients in a generalized sense. In particular,
our result covers the classical DiPerna-Lions flows and, we also obtain
the well-posedness for degenerate Fokker-Planck equations with irregular coefficients.
Moreover, a large deviation principle of Freidlin-Wenzell type for this type of SDEs is established.
\end{abstract}
\maketitle
\rm

\section{Introduction}

The celebrated DiPerna-Lions theory \cite{Di-Li} says that if a vector field $b\in W^{1,1}_{loc}(\mR^d)$
has bounded divergence and $\frac{b(x)}{1+|x|}\in L^1(\mR^d)+L^\infty(\mR^d)$,
then there exists a unique regular Lagrangian flow to  ordinary differential equation (ODE) in $\mR^d$:
\begin{equation}
\dif X_t(x)=b(X_t(x))\dif t,\ \ X_0(x)=x.\label{ODE}
\end{equation}
This theory was later extended to the case of  BV vector field by Ambrosio \cite{Am}.
Their methods were based on the connection between ODEs and transport or continuity equations.
Recently, Crippa and De Lellis \cite{Cr-De-Le} developed a more direct argument to treat this problem
by using the  Hardy-Littlewood maximal functions if $b$ is assumed to be in $W^{1,p}_{loc}(\mR^d)$ for some $p>1$.
Moreover, Cipriano and Cruzeiro \cite{Ci-Cr} studied the non-smooth flows associated to (\ref{ODE})
when the exponential of divergence of $b$ satisfies some $L^p(\mR^d,\mu)$-type hypothesis, where $\mu$ is the standard Gaussian
measure on $\mR^d$. Such theory has also been extended to the classical
Wiener space by Ambrosio and Figalli \cite{Am-Fi} (see also Fang and Luo \cite{Fa-Lu}).

We now turn to the following It\^o's stochastic differential equation (SDE) in $\mR^d$:
\begin{equation} \label{SDE}
 \dif X_t(x)=b(X_t(x))\dif t+\sigma(X_t(x))\dif W_t,\ \ X_0(x)=x,
\end{equation}
where $b:\mR^d\to\mR^d$ and $\sigma: \mR^d\to\mR^d\times\mR^m$ are two measurable functions,
and $(W_t)_{t\in[0,1]}$ is an $m$-dimensional standard Brownian motion on the classical Wiener space $(\Omega,\sF, P)$,
i.e., $\Omega$ is the space of all $\mR^m$-valued continuous functions on $[0,1]$, $\sF$ is the associated Borel $\sigma$-field,
$P$ is the standard Wiener measure. For a generic point $\omega\in\Omega$, $W_t(\omega)=\omega_t$ is the coordinate process.
Let $\sF_t$ be the natural Brownian filtration generated by $\{W_s, s\leq t\}$.

In \cite{Fi}, Figalli proved the well-posedness of martingale solutions for SDE (\ref{SDE})
with Sobolev coefficients by studying the associated Fokker-Planck equations. His strategy
is similar to \cite{Am}. Recently, we \cite{Zh2} gave a direct construction for the almost everywhere
stochastic flow of (\ref{SDE}) by using the same argument as in Crippa and De Lellis \cite{Cr-De-Le}.
Further more, through linearizing Brownian motion, we \cite{RenJie-Zh} also proved a classical limit theorem
that the solutions of ODE (\ref{ODE}) converges to the solutions of Stratonovich's SDEs in a generalized sense.
In the works of \cite{Cr-De-Le, Zh2, RenJie-Zh}, the vector field $b$ needs to be in $W^{1,q}_{loc}(\mR^d)$ for some $q>1$. In the non-degenerate and regular case of diffusion coefficients, there have been numerous results about the existence and uniqueness of strong solutions to SDE (\ref{SDE}) with singular drift $b$
(cf. \cite{Zv, Gy-Ma, Kr-Ro, Zh3}, etc.).

The present work is a continuation of \cite{Zh2} and \cite{RenJie-Zh}, and the main aims of this paper are two folds:
First, we try to relax the assumptions on the diffusion and drift coefficients so that
the diffusion coefficients can be discontinuous for Stratonovich SDEs, $b$ can be in $W^{1,1}_{loc}(\mR^d)$,
and the divergence of $b$ can be polynomial
growth. Secondly, we prove a Freidlin-Wentzell large deviation principle for SDEs with Sobolev coefficents.

In order to obtain a Freidlin-Wentzell large deviation estimate for SDE (\ref{SDE}) with discontinuous coefficients,
we shall employ the weak convergence method of Dupuis and Ellis \cite{Du-El}. This method has been proved to be very effective
for various stochastic systems (cf. \cite{Bo-Du-El, Bu-Du-Ma, Re-Xu-Zh}, etc.), where the key point is to use the variational
representation of certain exponential Brownian functionals (cf. \cite{Bo-Du, Zh5}) to prove an equivalent Laplace principle.

This paper is organized as follows: In Section 2, we state our main results.
In Section 3, some preliminaries are given. In Section 4, the well-posedness theorems are proven. In Section 5,
we shall prove a large deviation principle for SDE (\ref{SDE}).

\section{Statement of Main Results}

Let $\sM(\mR^d)$ be the total of all locally finite Borel measures on $\mR^d$.
For $p\geq 1$ and $\mu\in\sM(\mR^d)$, let $L^p_\mu=L^p_\mu(\mR^d)$ be the usual $L^p$-space over $(\mR^d,\mu)$
and $W^{p,k}_{loc}(\mR^d)$ the usual local Sobolev space.
If $\mu=\sL(\dif x)$ is the Lebesgue measure, we simply write $L^p_\mu=:L^p$.
For $R>0$, by $B_R$ we denote the ball in $\mR^d$ with center zero and radius $R$.

First of all, we introduce the following general notion about $\mu$-almost everywhere stochastic flow of
SDE (\ref{SDE}) (cf. \cite{Le-Li, Zh2}):

\bd\label{Def1}
Let $X_t(\omega,x)$ be a $\mR^d$-valued measurable
stochastic field on $[0,1]\times\Omega\times\mR^d$. For $\mu\in\sM(\mR^d)$,
we say $X$ a {\bf $\mu$-almost everywhere stochastic flow} of SDE (\ref{SDE}) corresponding to $(b,\sigma)$ if
\begin{enumerate}[{\bf(A)}]
\item for some $p\geq 1$, there exists a constant $K_{p}>0$
such that for any nonnegative measurable function $\varphi\in L^p_\mu(\mR^d)$,
\be\label{Den}
\sup_{t\in[0,1]}\mE\int_{\mR^d}\varphi(X_t(x))\mu(\dif x)\leq K_p\|\varphi\|_{L^p_\mu};
\ee
\item for  $\mu$-almost all $x\in\mR^d$, $t\mapsto X_t(x)$ is a continuous
($\sF_t$)-adapted  process satisfying  that
$$
\int^1_0|b(X_s(x))|\dif s+\int^1_0|\sigma(X_s(x))|^2\dif s<+\infty, \ \ P-a.s.,
$$
and
$$
X_t(x)=x+\int^t_0b(X_s(x))\dif s+\int^t_0\sigma(X_s(x))
\dif W_s,\ \ \forall t\in[0,1].
$$

\end{enumerate}
\ed

We first consider the following Stratonovich SDE:
\begin{equation*}
 \dif X_t(x)=b(X_t(x))\dif t+\sigma(X_t(x))\circ\dif W_t,\ \ X_0(x)=x,
\end{equation*}
or equivalent It\^o's form:
\begin{equation*}
 \dif X_t(x)=[b+\tfrac{1}{2}\sigma^{jl}\p_j\sigma^{\cdot l}](X_t(x))\dif t+\sigma(X_t(x))\dif W_t,\ \ X_0(x)=x.
\end{equation*}
Here and below, we use the conventions that the repeated indices in a product will be summed automatically,
and all derivatives and divergence are taken in the distributional sense.
By definitions, $\div\sigma^{\cdot l}:=\p_i\sigma^{il}$, $l=1,\cdots, m$

The following result is an extension of \cite[Theorem 2.6]{Zh2} to Stratonovich's SDE.
\bt\label{Main2}
Assume that for some $r\in[0,+\infty)$,
\begin{equation}
\frac{|b|+|\nabla\sigma|}{1+|x|}, |\sigma|\in L^\infty(B^c_r),\ \
b\in W^{1,1}_{loc}(\mR^d),\ \ \sigma\in W^{2,2}_{loc}(\mR^d), \label{Con2}
\end{equation}
and for some $\eps\in(0,1)$,
\begin{align}
 [\div b]^-,\ \ |\div\sigma|,\ \
\sup_{|z|\leq\eps}|\sigma(\cdot-z)|\cdot|\nabla\div\sigma|\in L^\infty(\mR^d).\label{Con1}
\end{align}
Then there exists a unique $\sL$-almost everywhere stochastic flow $X_t(x)$
in the sense of Definition \ref{Def1} corresponding to $(b_\sigma,\sigma)$ with $p=1$ in (\ref{Den}), where
$b_\sigma=b+\frac{1}{2}\sigma^{jl}\p_j\sigma^{\cdot l}$.
\et
\br
If $\div\sigma=\div b=0$, then from the proof below, one can see that
$$
\int_{\mR^d}\varphi(X_t(x))\dif x=\int_{\mR^d}\varphi(x)\dif x\ a.s.,\ \ \forall t\in[0,1],
$$
which means that stochastic flow $x\mapsto X_t(x)$ is incompressible. In this case, $b$ and $\sigma$
in Theorem \ref{Main2} only need to satisfy (\ref{Con2}) and so, are allowed to be singular in a finite ball.
If $\sigma$ vanishes, then our result covers the classical DiPerna-Lions flow.
\er

Our next aim is to relax the assumption $[\div b]^-\in L^\infty(\mR^d)$ so that $[\div b]^-$ can be polynomial growth.
We shall prove that:
\bt\label{Main1}
Assume that for some $q>1$,
\begin{equation}
|\nabla b|, |\nabla \sigma|^2\in L^q_{loc}(\mR^d), \ \
\frac{|b|+|\sigma|}{1+|x|}\in L^\infty(\mR^d),\label{Ep00}
\end{equation}
and there exist functions $\lambda\in C^2(\mR^d)$
and $\gamma_1,\gamma_2,\gamma_3$ satisfying that
for all small $y$ in $B_\eps$ and all $x\in\mR^d$,
\begin{align}
\lambda(x)\leq \gamma_1(x-y),\ \ |\nabla\lambda(x)|\leq\gamma_2(x-y),\ \ |\nabla^2\lambda(x)|\leq\gamma_3(x-y),\label{Ep0}
\end{align}
such that for all $p\geq1$,
\begin{align}
\int_{\mR^d}\exp\Big\{p\Big([\div b]^-+|b|\gamma_2+|\sigma|^2(\gamma^2_2+\gamma_3)
+|\nabla\sigma|^2\Big)(x)+\gamma_1(x)\Big\}\dif x<+\infty.\label{Ep1}
\end{align}
Let $\mu(\dif x)=e^{\lambda(x)}\dif x$. Then there exists a unique $\mu$-almost everywhere stochastic flow $X_t(x)$
in the sense of Definition \ref{Def1} corresponding to $(b,\sigma)$ with any $p>1$ in (\ref{Den}).
\et

\br
In this theorem, assumptions (\ref{Ep0}) and (\ref{Ep1}) are a little bit complicated.
We now explain them by introducing two examples.
\begin{enumerate}[(1)]
\item Let $\lambda(x)=-\a\log(1+|x|^2)$ for some $\a>\frac{d}{2}$.
For all $|y|\leq\frac{1}{2}$ and $x\in\mR^d$, we have
\begin{align*}
\lambda(x)&\leq-\a\log\left(1+(|x-y|-|y|)^2\right)\leq-\a\log\left(1+\tfrac{1}{2}|x-y|^2-|y|^2\right)\\
&\leq-\a\log\left(\tfrac{3}{4}+\tfrac{1}{2}|x-y|^2\right)
\leq-\a\log\left(1+|x-y|^2\right)+\a\log 2=:\gamma_1(x-y),
\end{align*}
and
$$
|\nabla\lambda(x)|\leq\frac{2\a|x|}{1+|x|^2}\leq\frac{4\a}{1+|x|}\leq\frac{8\a}{1+|x-y|}=:\gamma_2(x-y),
$$
$$
|\nabla^2\lambda(x)|\leq\frac{6\a}{1+|x|^2}\leq\frac{6\a}{1+\tfrac{1}{2}|x-y|^2-|y|^2}\leq\frac{12\a}{1+|x-y|^2}=:\gamma_3(x-y).
$$
In this case, if $b$ and $\sigma$ are linear growth, then condition (\ref{Ep1}) reduces to
$$
\int_{\mR^d}\frac{\exp\left\{p([\div b]^-+|\nabla\sigma|^2)(x)\right\}}{(1+|x|^2)^\a}\dif x<+\infty,\ \ \forall p\geq 1.
$$

\item Let $\lambda(x)=-|x|^{2\a}$ for some $\a\geq 1$. For all $|y|\leq\frac{1}{2}$ and $x\in\mR^d$, we have
$$
\lambda(x)\leq-(|x-y|-|y|)^{2\a}\leq-(|x-y|-\tfrac{1}{2})^{2\a}\leq C_\a-\tfrac{1}{2}|x-y|^{2\a}=:\gamma_1(x-y),
$$
and
$$
|\nabla\lambda(x)|\leq 2\a|x|^{2\a-1}\leq 2\a(|x-y|+\tfrac{1}{2})^{2\a-1}=:\gamma_2(x-y),
$$
$$
|\nabla^2\lambda(x)|\leq 4\a^2|x|^{2\a-2}\leq 4\a^2(|x-y|+\tfrac{1}{2})^{2\a-2}=:\gamma_3(x-y).
$$
In this case, if for some $\beta\in[0,1)$,
$$
\frac{|b(x)|}{1+|x|^\beta},\ \frac{|\sigma(x)|}{(1+|x|)^{\beta-\a}}\in L^\infty(\mR^d),
$$
then by Young's inequality, condition (\ref{Ep1}) reduces to
$$
\int_{\mR^d}\exp\left\{p([\div b]^-+|\nabla\sigma|^2)(x)-\tfrac{1}{4}|x|^{2\a}\right\}\dif x<+\infty, \ \ \forall p\geq 1.
$$
\end{enumerate}
\er
\br
Recently, Fang-Luo-Thalmaier \cite{Fa-Lu-Th} also studied the stochastic differential equations in Gaussian space
with Sobolev coefficients. However, our result is more general than
\cite[Theorem 1.3]{Fa-Lu-Th}. In particular, from Remark (1) above,
one can see that the condition 1.3 in \cite[Theorem 1.2]{Fa-Lu-Th} is not necessary.
\er
As an easy consequence of Theorem \ref{Main1} and \cite[Theorem 1.1]{Ro-Zh}, we have:
\bc\label{Co}
Assume that $b$ and $\sigma$ are bounded measurable functions and for some $q>1$,
$$
|\nabla b|, |\nabla \sigma|^2\in L^q_{loc}(\mR^d),
$$
and  (\ref{Ep1}) holds. Then for any probability density function $\phi$ with $\int_{\mR^d}\phi(x)^re^{(1-r)\lambda(x)}\dif x<+\infty$,
where $r>\frac{q}{q-1}=:p$, and $\lambda(x)$ is from Theorem \ref{Main1},
there exists a unique distribution solution to the following Fokker-Planck equation
\begin{align}
\p_t u_t=-\div(b u_t)+\tfrac{1}{2}\p^2_{ij}([\sigma^{il}\sigma^{jl}]u_t),\ \ u_0=\phi,\label{Eq1}
\end{align}
in the class of
$$
\cM_p:=\left\{u_t\in L^p_{loc}(\mR^d): u_t(x)\geq 0, \int_{\mR^d} u_t(x)\dif x=1,
\sup_{t\in[0,1]}\int_{\mR^d}u_t(x)^{p}e^{(1-p)\lambda(x)}\dif x<+\infty\right\}.
$$
\ec
\begin{proof}
Let $X_0$ be an $\sF_0$-measurable random variable with distribution $\phi(x)\dif x$. It is easy to see that
$Y_t:=X_t(X_0)$ solves the SDE:
$$
Y_t=X_0+\int^t_0 b(Y_s)\dif s+\int^t_0\sigma(Y_s)\dif W_s.
$$
Let $\mu(\dif x)=e^{\lambda(x)}\dif x$. Now for any $\varphi\in C^\infty_c(\mR^d)$, by H\"older's inequality, we have
\begin{align*}
\mE\varphi(Y_t)&=\mE(\mE\varphi(X_t(x))|x=X_0)=\int_{\mR^d}\mE\varphi(X_t(x)) \phi(x)\dif x\\
&\leq\left(\int_{\mR^d}|\mE\varphi(X_t(x))|^{\frac{r}{r-1}}\mu(\dif x)\right)^{1-\frac{1}{r}}\left(\int_{\mR^d}
\Big(\phi(x)e^{-r\lambda(x)}\mu(\dif x)\right)^{\frac{1}{r}}\\
&\leq\left(\mE\int_{\mR^d}|\varphi(X_t(x))|^{\frac{r}{r-1}}\mu(\dif x)\right)^{1-\frac{1}{r}}\left(\int_{\mR^d}
\phi(x)^re^{(1-r)\lambda(x)}\dif x\right)^{\frac{1}{r}}\leq C_\phi\|\varphi\|_{L^q_\mu}.
\end{align*}
Hence, there exists a $u\in\cM_p$ such that for any $\varphi\in C^\infty_c(\mR^d)$ and $t\in[0,1]$,
$$
\int_{\mR^d}\varphi(x) u_t(x)\dif x=\mE\varphi(Y_t)\leq C_\phi\|\varphi\|_{L^q_\mu}.
$$
By It\^o's formula, it is easy to check that $u$ is a distribution solution of (\ref{Eq1}).
The uniqueness follows from \cite[Theorem 1.1]{Ro-Zh}.
\end{proof}
\br
Compared with the result of Le Bris and Lions \cite[Proposition 5]{Le-Li1}, their well-posedness for equation (\ref{Eq1})
was given in the following space
$$
\{u\in L^\infty(0,1;(L^1\cap L^\infty)(\mR^d)), \sigma^\mathrm{t}\nabla u\in L^2(0,1; L^2(\mR^d))\}.
$$
Moreover, the conditions on $b$ and $\sigma$ are different.
\er
Next, we consider Freidlin-Wentzell's large deviation estimate of SDE (\ref{SDE}) in the situation of Theorem \ref{Main1}.
For $\eps\in(0,1)$, let $X_{\eps,t}(x)$ solve the following SDE in the sense of Definition \ref{Def1}:
\begin{equation}
\dif X_{\eps,t}(x)=b(X_{\eps,t}(x))\dif t+\sqrt{\eps}\sigma(X_{\eps,t}(x))\dif W_t,\ \ X_{\eps,0}(x)=x.\label{P5}
\end{equation}
We need to fix another weighted measure $\nu(\dif x)=e^{\rho(x)}\dif x$ such that
$$
\int_{\mR^d}|x|^{2p}\nu(\dif x)<+\infty,\ \ \forall p\geq 1.
$$
Thus we can consider equation (\ref{P5}) as an infinite-dimensional stochastic equation in Banach space $L^{2p}_\nu(\mR^d)$, $p\geq 1$:
$$
X_{\eps,t}=\mathrm{Id}+\int^t_0b(X_{\eps,s})\dif s+\sqrt{\eps}\int^t_0\sigma(X_{\eps,s})\dif W_s.
$$
The large deviation result is stated as follows.
\bt\label{Main3}
Assume that $b$ and $\sigma$ satisfy the same assumptions as in Theorem \ref{Main1}. Then
the family of random variables $(X_\eps)_{\eps\in(0,1)}$ as taking values in space
$\mS:=L^{2p}_\nu(\mR^d; C([0,1];\mR^d))$, $p\geq 1$ satisfies the large deviation principle.
More precisely, for any $B \in\cB(\mS)$, we have
$$
-\inf_{f\in B^o}I(f)\leq\varliminf_{\eps\rightarrow 0}\eps\log P(X_\eps\in B)
\leq\varlimsup_{\eps\rightarrow 0}\eps\log P(X_\eps\in B)\leq -\inf_{f\in \bar B}I(f),
$$
where $I(f):=\frac{1}{2}\inf_{\{h\in L^2(0,1):~f=X^h\}}\|h\|^2_{L^2}$, and $X^h$ solves the following equation:
\begin{equation}
X_t=\mathrm{Id}+\int^t_0b(X_{s})\dif s+\int^t_0\sigma(X_s) h_s\dif s.\label{C1}
\end{equation}
Here the closure and interior are taken in $\mS$.
\et
\br
Although Corollary \ref{Co} and Theorem \ref{Main3} are given under the assumptions of Theorem \ref{Main1},
similar results also hold for Stratonovich's SDE in the situation of Theorem \ref{Main2}.
\er
\section{Preliminaries}

\subsection{Two estimates on regular stochastic flows}
In this subsection, we assume that $b,\sigma\in C^\infty_b(\mR^d)$ are bounded and have bounded derivatives of all orders.
In this case, it is well known that SDE (\ref{SDE}) defines a
$C^\infty$-diffeomorphism flows $X_t(x), x\in\mR^d, t\in[0,1]$ (cf. \cite{Ik-Wa, Ku, Ma}).
We first recall the following well known result about the Jacobian determinant
(for example, see \cite[Lemma 3.1]{Zh2}).
\bl\label{Le5}
For any $t\in[0,1]$ and $x\in\mR^d$, we have
\begin{align}
\det(\nabla X_t(x))=\exp\left\{\int^t_0\div \sigma(X_s(x))\dif W_s+\int^t_0\Big[\div b
-\tfrac{1}{2}\p_i\sigma^{jl}\p_j\sigma^{il}\Big](X_s(x))\dif s\right\},\label{Es44}
\end{align}
and for any $p\geq 1$,
\begin{align}
\mE|\det(\nabla X^{-1}_t(x))|^p\leq\exp\left\{tp\Big(\|[-\div b+\tfrac{1}{2}\p_i\sigma^{jl}\p_j\sigma^{il}
+\sigma^{il}\p_{ij}^2\sigma^{jl}+\tfrac{p}{2}|\div \sigma|^2]^+\|_\infty\Big)\right\}.\label{Es5}
\end{align}
\el
Below, let $\lambda$ be a $C^2$-function on $\mR^d$ and define
\begin{align*}
\mu(\dif x):=e^{\lambda(x)}\dif x.
\end{align*}
We write
$$
\cJ_t(\omega,x):=\frac{(X_t(\omega,\cdot))_\sharp\mu(\dif x)}{\mu(\dif x)},\ \
\cJ^-_t(\omega,x):=\frac{(X^{-1}_t(\omega,\cdot))_\sharp\mu(\dif x)}{\mu(\dif x)},
$$
which means that for any nonnegative measurable function $\varphi$ on $\mR^d$,
\begin{align}
\int_{\mR^d}\varphi(X_t(\omega,x))\mu(\dif x)&=\int_{\mR^d}\varphi(x) \cJ_t(\omega,x)\mu(\dif x),\label{P0}\\
\int_{\mR^d}\varphi(X^{-1}_t(\omega, x))\mu(\dif x)&=\int_{\mR^d}\varphi(x) \cJ^-_t(\omega, x)\mu(\dif x).\label{P1}
\end{align}
It is easy to see that for almost all $\omega$ and all $(t,x)\in[0,1]\times\mR^d$,
\begin{align}
\cJ_t(\omega,x)=[\cJ^-_t(\omega,X^{-1}_t(\omega,x))]^{-1},\label{P2}
\end{align}
and by It\^o's formula and (\ref{Es44}),
\begin{align}
\cJ^-_t(x)=e^{\lambda(X_t(x))-\lambda(x)}\det(\nabla X_t(x))
=\exp\left\{\int^t_0\Lambda^{\sigma}_1(X_s(x))\dif W_s+\int^t_0\Lambda^{b,\sigma}_2(X_s(x))\dif s\right\},\label{P3}
\end{align}
where $\Lambda^\sigma_1(x):=\Big[\div\sigma+\sigma^{i\cdot}\p_i\lambda\Big](x)$ and
\begin{align*}
\Lambda^{b,\sigma}_2(x):=\Big[\div b+b^i\p_i\lambda+
\tfrac{1}{2}(\sigma^{il}\sigma^{jl}\p^2_{ij}
\lambda-\p_i\sigma^{jl}\p_j\sigma^{il})\Big](x).
\end{align*}

We now give an $L^p$-estimate for $\cJ_t(x)$, which is crucial for Theorem \ref{Main1} and inspired by \cite{Cr, Ci-Cr}.
\bl\label{Le00}
Assume that $\mu(\mR^d)<+\infty$. Then for any $t\in[0,1]$ and $p>1$, we have
\be
\mE\int_{\mR^d}|\cJ_t(x)|^p\mu(\dif x)\leq \mu(\mR^d)^{\frac{p}{p+1}}\Bigg(\sup_{t\in[0,1]}\int_{\mR^d}
\exp\Big\{tp^3|\Lambda^\sigma_1(x)|^2-tp^2 \Lambda^{b,\sigma}_2(x)
\Big\}\mu(\dif x)\Bigg)^{\frac{1}{p+1}}.\label{Es1}
\ee
\el
\begin{proof}
By (\ref{P1}) and (\ref{P2}), we have
\begin{align}
&\mE\int_{\mR^d}|\cJ_t(x)|^p\mu(\dif x)=\mE\int_{\mR^d}|\cJ^-_t(x)|^{1-p}\mu(\dif x).\label{Es2}
\end{align}
Since for any $\a\in\mR$,
$$
t\mapsto \exp\left\{\a\int^t_0\Lambda^\sigma_1(X_s(x))\dif W_s-\frac{\a^2}{2}
\int^t_0|\Lambda^\sigma_1(X_s(x))|^2\dif s\right\}
$$
is a continuous exponential martingale, by   (\ref{P3}) and H\"older's inequality,
for any $\a\in\mR$ and $q>1$, we have
$$
\mE|\cJ^-_t(x)|^{\a}\leq\left(\mE\exp\left\{\int^t_0\Big[\tfrac{q^2\a^2}{2(q-1)}|\Lambda^\sigma_1(X_s(x))|^2
+\a q\Lambda^{b,\sigma}_2(X_s(x))\Big]\dif s\right\}\right)^{\frac{1}{q}}.
$$
For the simplicity of notation, we write
\begin{align*}
\phi_{\a,q}(x):=\tfrac{q^2\a^2}{2(q-1)}|\Lambda^\sigma_1(x)|^2+\a q\Lambda^{b,\sigma}_2(x).
\end{align*}
By Jensen's inequality, we have
\begin{align*}
\mE\int_{\mR^d}|\cJ^-_t(x)|^{1-p}\mu(\dif x)&\leq
\int_{\mR^d}\left(\mE e^{\int^t_0\phi_{1-p,q}(X_s(x))\dif s}\right)^{\frac{1}{q}}\mu(\dif x)\\
&\leq\int_{\mR^d}\left(\frac{1}{t}\int^t_0\mE e^{t\phi_{1-p,q}(X_s(x))}\dif s\right)^{\frac{1}{q}}\mu(\dif x)\\
&\leq\mu(\mR^d)^{1-\frac{1}{q}}
\left(\frac{1}{t}\int^t_0\mE\int_{\mR^d} e^{t\phi_{1-p,q}(X_s(x))}\mu(\dif x)\dif s\right)^{\frac{1}{q}}\\
&\stackrel{(\ref{P0})}{=}\mu(\mR^d)^{1-\frac{1}{q}}\left(\frac{1}{t}\int^t_0
\mE\int_{\mR^d} e^{t\phi_{1-p,q}(x)}\cJ_s(x)
\mu(\dif x)\dif s\right)^{\frac{1}{q}}\\
&\leq\mu(\mR^d)^{1-\frac{1}{q}}\left(\int_{\mR^d} e^{\frac{pt}{p-1}\phi_{1-p,q}(x)}
\mu(\dif x)\right)^{\frac{p-1}{pq}}\\
&\quad\times\left[\sup_{s\in[0,1]}\mE\int_{\mR^d}|\cJ_s(x)|^p\mu(\dif x)\right]^{\frac{1}{pq}},
\end{align*}
which together with (\ref{Es2})  implies that
$$
\sup_{s\in[0,1]}\mE\int_{\mR^d}|\cJ_s(x)|^p\mu(\dif x)
\leq \mu(\mR^d)^{\frac{p(q-1)}{pq-1}}\left(\sup_{t\in[0,1]}\int_{\mR^d} e^{\frac{pt}{p-1}\phi_{1-p,q}(x)}\mu(\dif x)\right)^{\frac{p-1}{pq-1}}.
$$
The proof is complete by simplifying the above expression with $q=p$.
\end{proof}
\br
From (\ref{Es1}), one sees that by letting $p\downarrow 1$,
$$
\mE\int_{\mR^d}|\cJ_t(x)|\mu(\dif x)\leq \mu(\mR^d)^{\frac{1}{2}}\Bigg(\int_{\mR^d}
\exp\Big\{|\Lambda^\sigma_1(x)|^2+|\Lambda^{b,\sigma}_2(x)|\Big\}\mu(\dif x)\Bigg)^{\frac{1}{2}}.
$$
\er
\subsection{Two lemmas related to (\ref{Den})}

The following lemma will play a crucial role for taking limits below (cf. \cite{Zh2, RenJie-Zh}).
\bl\label{Le4}
Let $\mu\in\sM(\mR^d)$ and  $(X_n)_{n\in\mN}$ be a family of random fields on $\Omega\times\mR^d$.
Suppose that $X_n$ converges to $X$ for $P\otimes\mu$-almost all $(\omega,x)$, and for some $p\geq 1$,
there is a constant $K_p>0$ such that for any nonnegative
measurable function $\varphi\in L^p_\mu(\mR^d)$,
\begin{align}
\sup_n\mE\int_{\mR^d}\varphi(X_n(x))\mu(\dif x)\leq K_p\|\varphi\|_{L^p_\mu}.\label{Lp1}
\end{align}
Then we have:

(i). For any nonnegative measurable function $\varphi\in L^p_\mu(\mR^d)$,
\begin{align}
\mE\int_{\mR^d}\varphi(X(x))\mu(\dif x)\leq K_p\|\varphi\|_{L^p_\mu}.\label{Ps2}
\end{align}

(ii). If $\varphi_n$ converges to $\varphi$ in $L^p_\mu(\mR^d)$, then for any $N>0$,
\begin{align}
\lim_{n\to\infty}\mE\int_{B_N}|\varphi_n(X_n(x))-\varphi(X(x))|\mu(\dif x)=0.\label{Ps22}
\end{align}
\el
\begin{proof}
(i). First of all, for any nonnegative continuous function $\varphi\in C_c(\mR^d)$ with compact support,
by Fatou's lemma and (\ref{Lp1}), we have
$$
\mE\left(\int_{\mR^d}\varphi(X(x))\dif x\right)\leq
\varliminf_{n\to\infty}\mE\left(\int_{\mR^d}\varphi(X_n(x))\mu(\dif x)\right)\leq K_p\|\varphi\|_{L^p_\mu}.
$$
Let $O\subset\mR^d$ be a bounded open set. Define
$$
\varphi_n(x):=1-\left(\frac{1}{1+\mathrm{distance}(x,O^c)}\right)^n.
$$
Then $\varphi_n\in C_c(\mR^d)$ and for every $x\in\mR^d$,
$$
\varphi_n(x)\uparrow 1_O(x)\mbox{ as $n\to\infty$}.
$$
By the monotone convergence theorem, we find that (\ref{Ps2}) holds for $\varphi=1_O$.

We now extend (\ref{Ps2}) to the indicator function of any bounded Borel set. Without loss of generality, we consider
Borel sets in $(0,1]^d$, and define
$$
\sC:=\left\{A\in\cB((0,1]^d): \mE\left(\int_{\mR^d}1_A(X(x))\mu(\dif x)\right)\leq K_p \mu(A)^{1/p}\right\}
$$
and
$$
\sA:=\left\{A=\Pi_{i=1}^d(\a_i,\b_i]: 0<\a_i\leq \b_i\leq 1\right\}.
$$
It is easy to see that $\sC$ is a monotone class and $\sA$ is a semi-algebra on $(0,1]^d$.
Let $\sA_{\Sigma f}$ be the algebra generated by $\sA$ through finite disjoint unions.
Since all open subsets of $(0,1]^d$ belong to $\sC$, by another approximation, one finds that $\sA_{\Sigma f}\subset\sC$.
Hence, by the monotone class theorem,
$$
\cB((0,1]^d)\supset\sC\supset\sigma(\sA_{\Sigma f})=\cB((0,1]^d).
$$

Let $\varphi$ be a bounded nonnegative measurable function on some bounded open set $O$.
By Lusin's theorem, there exists a sequence of
bounded continuous functions $\varphi_\eps$ with  supports in $O$ such that
$$
\|\varphi_\eps\|_\infty\leq\|\varphi\|_\infty,\ \ \lim_{\eps\to 0}\mu(A_\eps)=0,
$$
where $A_\eps:=\{x\in \mR^d: \varphi(x)\not=\varphi_\eps(x)\}$.
Hence,
\begin{align*}
\mE\left(\int_{\mR^d}|\varphi-\varphi_\eps|(X(x))\mu(\dif x)\right)&\leq
2\|\varphi\|_\infty\mE\left(\int_{\mR^d}1_{A_\eps}(X(x))\mu(\dif x)\right)\\
&\leq 2\|\varphi\|_\infty K_p\mu(A_\eps)^{1/p}\stackrel{\eps\to 0}{\longrightarrow} 0.
\end{align*}
For general unbounded nonnegative measurable function $\varphi$ on $\mR^d$,
we can approximate it by the monotone convergence theorem again.

(ii). Let $\varphi_m\in C_c(\mR^d)$ converge to $\varphi$ in $L^p_\mu(\mR^d)$. By (\ref{Lp1}) and (\ref{Ps2}), we have
\begin{align*}
\mE\int_{B_N}|\varphi_n(X_n(x))-\varphi(X(x))|\mu(\dif x)
&\leq K_p\|\varphi_n-\varphi\|_{L^p_\mu}+\mE\int_{B_N}|\varphi(X_n(x))-\varphi(X(x))|\mu(\dif x)\\
&\leq K_p\|\varphi_n-\varphi\|_{L^p_\mu}+2K_p\|\varphi_m-\varphi\|_{L^p_\mu}\\
&\quad+\mE\int_{B_N}|\varphi_m(X_n(x))-\varphi_m(X(x))|\mu(\dif x),
\end{align*}
which converges to zero by first letting $n\to \infty$ and then $m\to\infty$.
\end{proof}

Let $\varrho\geq 0$ be a smooth function in $\mR^d$ with supp${\varrho}\subset B_1$ and $\int_{\mR^d}\varrho(x)\dif x=1$.
For $\eps>0$, set
\begin{align}
\varrho_\eps(x):=\eps^{-d}\varrho(\eps^{-1}x).\label{Mol}
\end{align}
For a function $b\in L^1_{loc}(\mR^d)$, define
\begin{align}
b_\eps(x):=b*\varrho_\eps(x)=\int_{\mR^d}b(y)\varrho_\eps(x-y)\dif y,\label{Es01}
\end{align}
and for any $R>0$ and $\varphi\in L^1_{loc}(\mR^d)$,
$$
M_R\varphi(x):=\sup_{0<s<R}\fint_{B_s}\varphi(x+y)\dif y,
$$
where $\fint_{B_s}\varphi(x+y)\dif y:=\frac{1}{|B_s|}\int_{B_s}\varphi(x+y)\dif y$.

We have the following elementary estimate.
\bl\label{Le01}
Let $b\in W^{1,1}_{loc}(\mR^d)$. Then there exists an $\sL$-null set $A\subset\mR^d$
such that for all $x,y\notin A$,
$$
|b(x)-b(y)|\leq 2^d\int^{|x-y|}_0\!\!\!\fint_{B_s}|\nabla b|(x+z) \dif z\dif s
+2^d\int^{|x-y|}_0\!\!\!\fint_{B_s}|\nabla b|(y+z) \dif z\dif s.
$$
In particular, for any $R>0$ and $x,y\notin A$ with $|x-y|\leq R$,
\begin{align}
|b(x)-b(y)|\leq 2^d|x-y|(M_R|\nabla b|(x)+(M_R|\nabla b|(y)).\label{Epp2}
\end{align}
\el
\begin{proof}  Let $b_\eps(x)$ be defined by (\ref{Es01}). For $r>0$, let
$\Pi(\dif z)$ denote the surface measure on the ball $\{z\in\mR^d: |z|=r\}$. Noting that
$$
|b_\eps(x)-b_\eps(x+z)|\leq |z|\int^1_0|\nabla b_\eps|(x+sz)\dif s,
$$
we have
\begin{align*}
\int_{|z|=r}|b_\eps(x)-b_\eps(x+z)|\Pi(\dif z)&\leq r\int^1_0\int_{|z|=r}|\nabla b_\eps|(x+sz)\Pi(\dif z)\dif s\\
&=r\int^1_0s^{1-d}\int_{|z|=sr}|\nabla b_\eps|(x+z) \Pi(\dif z)\dif s.
\end{align*}
Hence, for any $\ell>0$,
\begin{align*}
\int_{B_\ell}|b_\eps(x)-b_\eps(x+z)|\dif z&=
\int^\ell_0\int_{|z|=r}|b_\eps(x)-b_\eps(x+z)|\Pi(\dif z)\dif r\\
&\leq\int^\ell_0t\int^1_0s^{1-d}\int_{|z|=sr}|\nabla b_\eps|(x+z) \Pi(\dif z)\dif s\dif r\\
&=\int^1_0s^{-1-d}\int^{s\ell}_0r\int_{|z|=r}|\nabla b_\eps|(x+z) \Pi(\dif z)\dif r\dif s\\
&\leq\int^1_0s^{-d}\ell \int_{B_{s\ell}}|\nabla b_\eps|(x+z) \dif z\dif s\\
&=\ell^d\int^\ell_0s^{-d} \int_{B_s}|\nabla b_\eps|(x+z) \dif z\dif s.
\end{align*}
For any $x,y\in\mR^d$, set $\ell:=|x-y|$, then
\begin{align}
|b_\eps(x)-b_\eps(y)|&\leq \fint_{B_{\ell/2}}|b_\eps(x)-b_\eps(\tfrac{x+y}{2}+z)|\dif z
+\fint_{B_{\ell/2}}|b_\eps(y)-b_\eps(\tfrac{x+y}{2}+z)|\dif z\no\\
&\leq 2^d\fint_{B_\ell}|b_\eps(x)-b_\eps(x+z)|\dif z+2^d\fint_{B_{\ell}}|b_\eps(y)-b_\eps(y+z)|\dif z\no\\
&\leq 2^d\int^\ell_0\!\!\!\fint_{B_s}|\nabla b_\eps|(x+z) \dif z\dif s
+2^d\int^{\ell}_0\!\!\!\fint_{B_s}|\nabla b_\eps|(y+z) \dif z\dif s.\label{Es}
\end{align}
Since for any $R,\ell>0$,
$$
\lim_{\eps\to 0}\int^1_0\!\!\!\int_{B_R}|b_\eps-b|(x)\dif x\dif t=0
$$
and
$$
\lim_{\eps\to 0}\int^1_0\!\!\!\int_{B_R}\left(\int^\ell_0\!\!\!\fint_{B_s}|\nabla(b_\eps-b)|(x+z) \dif z\dif s\right)
\dif x\dif t=0,
$$
we can take limits $\eps\to 0$ for (\ref{Es}) and obtain the desired estimate.
\end{proof}

\bl\label{Le20}
Let $b\in W^{1,1}_{loc}(\mR^d)$. There exists an $\sL$-null set $A\subset\mR^d$
such that for any $\delta,\eps\in(0,\frac{1}{4})$, and all $x,y\in \mR^d\setminus A$
with $|x-y|\leq\sqrt{\delta}$,
\begin{align}
\frac{|b(x)-b(y)|}{\sqrt{|x-y|^2+\delta^2}}\leq 2^d(f_{\delta,\eps}(x)+f_{\delta,\eps}(y)),\label{Op4}
\end{align}
where
\begin{align*}
f_{\delta,\eps}(x)&:=\eps^{-d}\|\varrho\|_\infty\int_{B_1}|\nabla b|(x+z)\dif z+
\frac{1}{\delta}\int^\delta_0\!\!\!\fint_{B_s}|\nabla b|(x+z) \dif z\dif s\\
&\qquad\quad+\int^{\sqrt{\delta}}_\delta\frac{1}{s}\left(\fint_{B_s}|\nabla(b_\eps-b)|(x+z) \dif z\right)\dif s,
\end{align*}
and $b_\eps(x)=b*\varrho_\eps(x)$ is the mollifying vector field. Moreover, for any $R>0$,
\begin{align}
\int_{B_R}f_{\delta,\eps}(x)\dif x&\leq C_{\varrho,d}\eps^{-d}\|\nabla b\|_{L^1(B_{R+1})}
+\frac{\log \delta^{-1}}{2}\|\nabla(b_\eps-b)\|_{L^1(B_{R+1})},
\label{Op3}
\end{align}
where $C_{\varrho,d}$ only depends on $\|\varrho\|_\infty$ and $d$.
\el

\begin{proof}
Set $\ell:=|x-y|\leq\sqrt{\delta}$. By Lemma \ref{Le01}, we have
\begin{align*}
\frac{|b(x)-b(y)|}{\sqrt{|x-y|^2+\delta^2}}\leq
2^d\left(\frac{1}{\delta}\wedge\frac{1}{\ell}\right)
\left(\int^\ell_0\!\!\!\fint_{B_s}|\nabla b|(x+z) \dif z\dif s
+\int^{\ell}_0\!\!\!\fint_{B_s}|\nabla b|(y+z) \dif z\dif s\right).
\end{align*}
We make the following estimate:
\begin{align*}
\left(\frac{1}{\delta}\wedge\frac{1}{\ell}\right)\int^\ell_0\!\!\!\fint_{B_s}|\nabla b|(x+z) \dif z\dif s&\leq
\frac{1}{\delta}\int^\delta_0\!\!\!\fint_{B_s}|\nabla b|(x+z) \dif z\dif s
+\frac{1_{\ell>\delta}}{\ell}\int^\ell_\delta\!\!\!\fint_{B_s}|\nabla b|(x+z) \dif z\dif s\\
&\leq\frac{1}{\delta}\int^\delta_0\!\!\!\fint_{B_s}|\nabla b|(x+z) \dif z\dif s
+\frac{1}{\ell}\int^\ell_\delta\!\!\!\fint_{B_s}|\nabla b_\eps|(x+z) \dif z\dif s\\
&\quad+\frac{1_{\ell>\delta}}{\ell}\int^\ell_\delta\!\!\!\fint_{B_s}|\nabla(b_\eps-b)|(x+z) \dif z\dif s\\
&\leq\frac{1}{\delta}\int^\delta_0\!\!\!\fint_{B_s}|\nabla b|(x+z) \dif z\dif s
+\sup_{z\in B_{\sqrt{\delta}}}|\nabla b_\eps(x+z)|\\
&\quad+\int^{\sqrt{\delta}}_\delta\frac{1}{s}\left(\fint_{B_s}|\nabla(b_\eps-b)|(x+z) \dif z\right)\dif s.
\end{align*}
Estimate (\ref{Op4}) now follows by noting that
$$
\sup_{z\in B_{\sqrt{\delta}}}|\nabla b_\eps|(x+z)\leq\eps^{-d}\|\varrho\|_\infty\int_{B_1}|\nabla b|(x+z) \dif z
$$
provided that $\eps,\delta<\frac{1}{4}$.

As for (\ref{Op3}), by Fubini's theorem, we have
\begin{align*}
\int^1_0\!\!\!\int_{B_R}f_{\delta,\eps}(x)\dif x\dif s
&\leq \eps^{-d}\|\varrho\|_\infty\int_{B_R}\!\int_{B_1}|\nabla b|(x+z) \dif z\dif x
+\int^1_0\!\!\!\int_{B_{R+1}}|\nabla b|(z) \dif z\dif t\\
&\quad+\int^{\sqrt{\delta}}_\delta\frac{1}{s}\dif s\int^1_0\!\!\!\int_{B_{R+1}}|\nabla(b_\eps-b)|(z) \dif z\dif t\\
&\leq(\eps^{-d}\|\varrho\|_\infty|B_1|+1)\int^1_0\!\!\!\int_{B_{R+1}}|\nabla b|(z) \dif z\dif t\\
&\quad+\log\left(\frac{1}{\sqrt{\delta}}\right)\int^1_0\!\!\!\int_{B_{R+1}}|\nabla(b_\eps-b)|(z) \dif z\dif t.
\end{align*}
The proof is complete.
\end{proof}
We also recall the following well known result (cf. \cite{St}).
\bl\label{Le0}
For any $p>1$, there exists $C_{d,p}>0$ such that for any $N,R>0$ and $\varphi\in L^p_{loc}(\mR^d)$,
\begin{align}
\int_{B_N}(M_R\varphi(x))^p\dif x\leq{C_{d,p}}\int_{B_{N+R}}|\varphi(x)|^p\dif x.\label{Epp1}
\end{align}
\el

\subsection{An abstract criterion for Laplace principle}
Let $\mH$ be the Cameron-Martin space over the classical Wiener space, the space of all absolutely continuous functions
from $[0,1]$ to $\mR^d$, which is isomorphic to $L^2(0,1;\mR^d)$
through the mapping $h\mapsto \int^\cdot_0 h_s\dif s$. Below, we always regard $\mH$ as $L^2(0,1;\mR^d)$.
For $M>0$,  set
$$
\cD_M:=\{h\in \mH: \|h\|_\mH\leq M\}
$$
and
\be
\cA_M:=\left\{
\begin{aligned}
&\mbox{ $h: [0,1]\to \mH$ is a simple and $(\sF_t)$-adapted }\\
&\mbox{ process, and for almost all $\omega$},\ \ h(\cdot,\omega)\in\cD_M
\end{aligned}
\right\}.\label{Op2}
\ee
We equip $\cD_M$ with the  weak convergence topology in $\mH$ so that
$\cD_M$ becomes a compact Polish space.
Let $\mS$ be a Polish space. A function $I: \mS\to[0,\infty]$ is given.
\bd
The function $I$  is called a rate function if for every $a<\infty$, the set
$\{f\in\mS: I(f)\leq a\}$ is compact in $\mS$.
\ed

Let $\{Z^\eps: \Omega\to\mS,\eps\in(0,1)\}$
be a family of measurable mappings. Assume that
there is a measurable map $Z_0: \mH\to \mS$ such that
\begin{enumerate}[{\bf (LD)$_\mathbf{1}$}]
\item
For any $M>0$, if a family $\{h_\eps, \eps\in(0,1)\}\subset\cA_M$ (as random variables in $\cD_M$)
converges in distribution to $h\in \cA_M$, then
for some subsequence $\eps_k$, $Z^{\eps_k}\Big(\cdot+\frac{1}{\sqrt{\eps_k}}\int^\cdot_0h^{\eps_k}_s(\cdot)\dif s\Big)$
converges in distribution to $Z_0(h)$ in $\mS$.
\end{enumerate}
\begin{enumerate}[{\bf (LD)$_\mathbf{2}$}]
\item
 For any $M>0$, if $\{h_n,n\in\mN\}\subset \cD_M$ weakly converges to $h\in\mH$,
then for some subsequence $h_{n_k}$, $Z_0(h_{n_k})$ converges to $Z_0(h)$ in $\mS$.
\end{enumerate}

For each $f\in\mS$, define
\be
I(f):=\tfrac{1}{2}\inf_{\{h\in\mH:~f=Z_0(h)\}}\|h\|^2_{\mH},\label{ra}
\ee
where $\inf\emptyset=\infty$ by convention. Then under {\bf (LD)$_\mathbf{2}$},
$I(f)$ is a  rate function.

We recall the following result due to \cite{bd0} (see also \cite[Theorem 4.4]{Zh5}).
\bt\label{Th2}
Under {\bf (LD)$_\mathbf{1}$} and {\bf (LD)$_\mathbf{2}$},
$\{Z^\eps,\eps\in(0,1)\}$ satisfies the Laplace principle with
the rate function $I(f)$ given by (\ref{ra}). More precisely, for each real bounded continuous
function  $g$ on $\mS$:
\be
\lim_{\eps\rightarrow 0}\eps\log\mE\left(\exp\left[-\frac{g(Z^\eps)}{\eps}\right]\right)
=-\inf_{f\in\mS}\{g(f)+I(f)\}.\label{La}
\ee
In particular, the family $\{Z^\eps,\eps\in(0,1)\}$
satisfies  the large deviation principle in $(\mS,\cB(\mS))$ with the rate function $I(f)$.
\et

\section{Proofs of Theorems \ref{Main2} and \ref{Main1}}

We first establish the following key stability estimate.
\bl\label{Le1}
Assume that for some $q\geq 1$,
$$
b, \hat b\in L^q_{loc}(\mR^d),\ \ |\nabla b|\in L^q_{loc}(\mR^d)
$$
and
$$
\sigma,\hat \sigma\in L^{2q}_{loc}(\mR^d),\ \ |\nabla\sigma|\in L^{2q}_{loc}(\mR^d).
$$
Let $\mu(\dif x)=e^{\lambda(x)}\dif x$ with $\lambda\in C(\mR^d)$.
Let $X_t(x)$ and $\hat X_t(x)$ be two $\mu$-almost everywhere stochastic flows of (\ref{SDE})
corresponding to $(b,\sigma)$ and $(\hat b,\hat \sigma)$ in the sense of Definition
\ref{Def1} with $p=q$ in (\ref{Den}). Then for any $N,R>1$ and $\eta, \delta,\eps\in(0,1)$,
there exist constants $C_1, C_2, C_3>0$ such that
\begin{align*}
\mE\int_{B_N}\left(\sup_{t\in[0,1]}|X_t(x)-\hat X_t(x)|^{2}\wedge 1\right)\mu(\dif x)&\leq \eta
+\frac{2\mu(B_N)}{R\eta}\mE\int_{B_N}\left(\sup_{t\in[0,1]}|X_t(x)|\vee|\hat X_t(x)|\right)\mu(\dif x)\\
&\quad+\frac{C_1(\eps^{-d}1_{q=1}+1_{q>1})}{\eta\log\delta^{-1}}
+\frac{C_2}{\eta}\|\nabla(b_\eps-b)\|_{L^1(B_{R+1})}1_{q=1}\\
&\quad+\frac{C_3}{\eta\delta\log\delta^{-1}}\left(\|b-\hat b\|_{L^q(B_{R})}
+\|\sigma-\hat\sigma\|_{L^{2q}(B_{R})}\right),
\end{align*}
where  $b_\eps(x)=b*\varrho_\eps(x)$, $C_1=C(R,N,\|\nabla b\|_{L^q(B_{R+1})},
\|\nabla\sigma\|_{L^{2q}(B_{R+1})}, K_q,\lambda)$ and $C_2=C_3=C(R,N, K_q,\lambda)$.
Here, $K_q$ is from (\ref{Den}).
\el
\begin{proof}
For $\delta>0$, let $\xi_\delta:\mR_+\to\mR_+$
be a smooth function with $0\leq \xi_\delta'(s)\leq 1$, $0\leq\xi''_\delta(s)\leq\frac{4}{\delta}$ and
$$
\xi_\delta(s)=
\left\{
\begin{aligned}
&s,&\ \ \  s\in[0,\delta/4];\\
&\delta/2,&\ \ \  s\in[\delta,\infty).
\end{aligned}
\right.
$$
By elementary calculations, we have
\begin{align}
s\leq 2\xi_\delta(s),\ \ s\in[0,\delta].\label{Ep2}
\end{align}
Set
$$
Z_t(\omega,x):=X_t(\omega,x)-\hat X_t(\omega, x)
$$
and
$$
\Phi(\omega,x):=\sup_{t\in[0,1]}\xi_\delta(|Z_t(\omega,x)|^2).
$$
We divide the proof into two steps.

{\bf (Step 1)}. In this step we prove that for any $N,R>1$, there exist constants $C_1,C_2, C_3>0$
as in the statement of the theorem such that for all $\delta,\eps\in(0,1)$,
\begin{align}
\mE\int_{B_N\cap G_R}
\log\left(\frac{\Phi(x)}{\delta^2}+1\right)\mu(\dif x)&\leq C_1\eps^{-d}+C_2\log\delta^{-1}
\int_{B_{R+1}}|\nabla(b_\eps-b)|(z)\dif z\no\\
&\quad+\frac{C_3}{\delta}\left(\|b-\hat b\|_{L^q(B_{R})}+\|\sigma-\hat\sigma\|_{L^{2q}(B_{R})}\right),\label{Es8}
\end{align}
where $G_R(\omega):=\Big\{x\in\mR^d: \sup_{t\in[0,1]}|X_t(\omega,x)|\vee|\hat X_t(\omega,x)|\leq R\Big\}$.

Noticing that for $\mu$-almost all $x\in\mR^d$ and all $t\in[0,1]$
$$
Z_t(x)=\int^t_0(b(X_s(x))-\hat b(\hat X_s(x)))\dif s+\int^t_0(\sigma(X_s(x))-\hat\sigma(\hat X_s(x)))\dif W_s,
$$
by It\^o's formula, we have
\begin{align*}
\log\left(\frac{\xi_\delta(|Z_t(x)|^2)}{\delta^2}+1\right)
&=2\int^t_0\frac{\xi'_\delta(|Z_s(x)|^2)\<Z_s(x), b(X_s(x))-\hat b(\hat X_s(x))\>}{\xi_\delta(|Z_s(x)|^2)+\delta^2}\dif s\\
&\quad+2\int^t_0\frac{\xi'_\delta(|Z_s(x)|^2)\<Z_s(x), (\sigma(X_s(x))-\hat \sigma(\hat X_s(x)))\dif W_s\>}
{\xi_\delta(|Z_s(x)|^2)+\delta^2}\\
&\quad+\int^t_0\frac{\xi'_\delta(|Z_s(x)|^2)\|\sigma(X_s(x))-\hat \sigma(\hat X_s(x))\|^2}
{\xi_\delta(|Z_s(x)|^2)+\delta^2}\dif s\\
&\quad+2\int^t_0\frac{\xi''_\delta(|Z_s(x)|^2)|(\sigma(X_s(x))
-\hat \sigma(\hat X_s(x)))^\mathrm{t}\cdot Z_s(x)|^2}{\xi_\delta(|Z_s(x)|^2)+\delta^2}\dif s\\
&\quad-2\int^t_0\frac{(\xi'_\delta(|Z_s(x)|^2))^2|(\sigma(X_s(x))
-\hat \sigma(\hat X_s(x)))^\mathrm{t}\cdot Z_s(x)|^2}{(\xi_\delta(|Z_s(x)|^2)+\delta^2)^2}\dif s\\
&=:I_1(t,x)+I_2(t,x)+I_3(t,x)+I_4(t,x)+I_5(t,x).
\end{align*}
Since $I_5(t,x)$ is negative, we can drop it. For $I_1(t,x)$, by (\ref{Ep2}), we have
\begin{align*}
\sup_{t\in[0,1]}|I_1(t,x)|&\leq4\int^1_0\frac{|b(X_s(x))-b(\hat X_s(x))|
\cdot 1_{|Z_s(x)|\leq\sqrt{\delta}}}{\sqrt{|Z_s(x)|^2+\delta^2}}\dif s\\
&\quad+\frac{2}{\delta}\int^1_0|b(\hat X_s(x))-\hat b(\hat X_s(x))|\dif s\\
&=:I_{11}(x)+I_{12}(x).
\end{align*}
Noting that
$$
G_R(\omega)\subset \{x: |X_t(\omega,x)|\leq R\}\cap\{x: |\hat X_t(\omega,x)|\leq R\},\ \ \forall t\in[0,1],
$$
by (\ref{Den}), we have
\begin{align}
\mE\int_{G_R}|I_{12}(x)|\mu(\dif x)&\leq
\frac{2}{\delta}\mE\int^1_0\!\!\!\int_{\mR^d}|1_{B_R}(b-\hat b)|(\hat X_s(x))\mu(\dif x)\dif s\no\\
&\leq\frac{2K_q}{\delta}\|1_{B_R}(b-\hat b)\|_{L^q_\mu}
\leq\frac{C_{q,R,\lambda}}{\delta}\|b-\hat b\|_{L^q(B_R)}.\label{PL0}
\end{align}
For $I_{11}(x)$, if $q=1$, by Lemma  \ref{Le20}, we have
\begin{align}
\mE\int_{G_R}|I_{11}(x)|\mu(\dif x)&\leq
2^{d+2}\mE\int^1_0\!\!\!\int_{G_R}[f_{\delta,\eps}(X_s(x))+f_{\delta,\eps}(\hat X_s(x))]\mu(\dif x)\dif s\no\\
&\leq C_d\int_{B_R}f_{\delta,\eps}(x)\mu(\dif x)\leq C_{d,R,\lambda}\int_{B_R}f_{\delta,\eps}(x)\dif x\no\\
&\leq C_{d,R,\lambda,\varrho}\left(\eps^{-d}\|\nabla b\|_{L^1(B_{R+1})}+
\log\delta^{-1}\|\nabla(b_\eps-b)\|_{L^1(B_{R+1})}\right);\label{PL00}
\end{align}
if $q>1$, by Lemma \ref{Le0}, we have
\begin{align}
\mE\int_{G_R}|I_{11}(x)|\mu(\dif x)&\leq C\mE\int^1_0\!\!\!\int_{G_R}(M_{\sqrt{\delta}}|\nabla b|(X_s(x))
+M_{\sqrt{\delta}}|\nabla b|(\hat X_s(x)))\mu(\dif x)\dif s\no\\
&\leq C\left(\int_{B_R}(M_{\sqrt{\delta}}|\nabla b|(x))^q\mu(\dif x)\right)^{1/q}\leq C \|\nabla b\|_{L^q(B_{R+1})}.\label{PL1}
\end{align}
For $I_2(t,x)$, set
$$
\tau_R(\omega,x):=\inf\Big\{t\in[0,1]: |X_t(\omega,x)|\vee \hat X_t(\omega,x)>R\Big\},
$$
then
$$
G_R(\omega)=\{x:\tau_R(\omega,x)=1\}.
$$
By Burkholder's inequality, Fubini's theorem and (\ref{Ep2}), we have
\begin{align*}
&\mE\int_{B_N\cap G_R}\sup_{t\in[0,1]}|I_{2}(t,x)|\mu(\dif x)\\
&\qquad\leq
\int_{B_N}\mE\left(\sup_{t\in[0,\tau_R(x)]}\left|\int^{t}_0\frac{\xi'_\delta(|Z_s(x)|^2)\<Z_s(x),(\sigma(X_s(x))
-\hat \sigma(\hat X_s(x)))\dif W_s\>}{\xi_\delta(|Z_s(x)|^2)+\delta^2}\right|\right)\mu(\dif x)\\
&\qquad\leq C\int_{B_N}\mE\left[\int^{ \tau_R(x)}_0\frac{(\xi'_\delta(|Z_s(x)|^2))^2
|Z_s(x)|^2|\sigma(X_s(x))-\hat \sigma(\hat X_s(x))|^2}
{(\xi_\delta(|Z_s(x)|^2)+\delta^2)^2}\dif s\right]^{\frac{1}{2}}\mu(\dif x)\\
&\qquad\leq C\mu(B_N)^{\frac{1}{2}}\left[\mE\int^1_0\!\!\!\int_{B_N\cap G_R}
\frac{|\sigma(X_s(x))-\hat \sigma(\hat X_s(x))|^2\cdot 1_{|Z_s(x)|\leq \sqrt{\delta}}}
{|Z_s(x)|^2+\delta^2}\mu(\dif x)\dif s\right]^{\frac{1}{2}}.
\end{align*}
As the treatment of $I_1(t,x)$, by Lemma \ref{Le0}, we can prove that
\begin{align}
\mE\int_{B_N\cap G_R}\sup_{t\in[0,1]}|I_{2}(t,x)|\mu(\dif x)
\leq C\|\nabla\sigma\|_{L^{2q}(B_{R+1})}
+\frac{C}{\delta}\|\sigma-\hat \sigma\|_{L^{2q}(B_{R})}\label{PL2}
\end{align}
and similarly,
\begin{align}
\mE\int_{B_N\cap G_R}\sup_{t\in[0,1]}|I_{3}(t,x)|\mu(\dif x)
&\leq C\|\nabla\sigma\|_{L^{2q}(B_{R+1})}
+\frac{C}{\delta}\|\sigma-\hat \sigma\|_{L^{2q}(B_{R})},\\
\mE\int_{B_N\cap G_R}\sup_{t\in[0,1]}|I_{4}(t,x)|\mu(\dif x)
&\leq C\|\nabla \sigma\|_{L^{2q}(B_{R+1})}
+\frac{C}{\delta}\|\sigma-\hat \sigma\|_{L^{2q}(B_{R})}.\label{PL3}
\end{align}
Combining (\ref{PL0})-(\ref{PL3}), we obtain (\ref{Es8}).

{\bf (Step 2)}. For any $\eta>0$, we have
\begin{align}
\mE\int_{B_N}\left(\Phi(x)\wedge 1\right)\mu(\dif x)
&\leq\eta+\mu(B_N)P\left\{\int_{B_N}\left(\Phi(x)\wedge 1\right)\mu(\dif x)\geq \eta\right\}\no\\
&\leq\eta+\mu(B_N)P\left\{\int_{B_N\cap G^c_R}\left(\Phi(x)\wedge 1\right)\mu(\dif x)\geq \frac{\eta}{2}\right\}\no\\
&\quad+\mu(B_N)P\left\{\int_{B_N\cap G_R}\left(\Phi(x)\wedge 1\right)\mu(\dif x)\geq \frac{\eta}{2}\right\}.\label{Es7}
\end{align}
By Chebyshev's inequality, we have
\begin{align}
P\left\{\int_{B_N\cap G^c_R}\left(\Phi(x)\wedge 1\right)\mu(\dif x)\geq \frac{\eta}{2}\right\}
&\leq P\left\{\mu(B_N\cap G^c_R)\geq \frac{\eta}{2}\right\}\leq \frac{2}{\eta}\mE\mu(B_N\cap G^c_R)\no\\
&\leq\frac{2}{R\eta}\mE\int_{B_N}\left(\sup_{t\in[0,1]}|X_t(x)|\vee|\hat X_t(x)|\right)\mu(\dif x).\label{Es9}
\end{align}
Set now
$$
\Psi_\delta(x):=\log\left(\frac{\Phi(x)}{\delta^2}+1\right).
$$
Notice that if $\Psi_\delta(x)\leq\frac{\log\delta^{-1}}{2}$, then $\Phi(x)<\delta$.
Hence, for any $\delta<\frac{\eta}{4\mu(B_N)}$, we have
\begin{align}
P\left\{\int_{B_N\cap G_R}\left(\Phi(x)\wedge 1\right)\mu(\dif x)\geq \frac{\eta}{2}\right\}
&\leq P\left\{\int_{B_N\cap G_R}\left(\Phi(x)\wedge 1\right)\cdot
1_{\{2\Psi_\delta(x)>\log{\delta^{-1}}\}}\mu(\dif x)\geq \frac{\eta}{4}\right\}\no\\
&+P\left\{\int_{B_N\cap G_R}\left(\Phi(x)\wedge 1\right)\cdot
1_{\{2\Psi_\delta(x)\leq \log{\delta^{-1}}\}}\mu(\dif x)\geq \frac{\eta}{4}\right\}\no\\
&\leq P\left\{\int_{B_N\cap G_R}\Psi_\delta(x)\mu(\dif x)\geq \frac{\eta\log{\delta^{-1}}}{8}\right\}+0\no\\
&\leq \frac{8}{\eta\log{\delta^{-1}}}\mE\int_{B_N\cap G_R}\Psi_\delta(x)\mu(\dif x).\label{Es6}
\end{align}
The result now follows by combining (\ref{Es8}), (\ref{Es7}), (\ref{Es9}) and (\ref{Es6}).
\end{proof}

Let $\chi\in C^\infty(\mR^d)$ be a nonnegative cutoff function
with
\begin{equation} \label{eq:4}
\|\chi\|_\infty\leq 1,\ \
\chi(x)=\left\{ \begin{aligned}
&1,\ \ |x|\leq 1, \\
&0,\ \ |x|\geq 2.
\end{aligned} \right.
\end{equation}
Set $\chi_n(x):=\chi(x/n)$ and define
\be
b_n:=b*\rho_n\cdot\chi_n,\ \ \sigma_n:=\sigma*\rho_n\cdot\chi_n,\label{Es3}
\ee
where $\rho_n=\varrho_{1/n}$ is the mollifiers given by (\ref{Mol}).

We are now in a position to give the proofs of Theorems \ref{Main2} and \ref{Main1}.

\begin{proof}[Proof of Theorem \ref{Main2}] Let $b_n$ and $\sigma_n$ be defined by (\ref{Es3}). Let $X^n_t(x)$ be the solution of
the Stratonovich SDE:
\begin{align*}
X^n_t(x)&=x+\int^t_0b_n(X^n_s(x))\dif s+\int^t_0\sigma_n(X^n_s(x))\circ\dif W_s\\
&=x+\int^t_0\tilde b_n(X^n_s(x))\dif s+\int^t_0\sigma_n(X^n_s(x))\dif W_s,
\end{align*}
where $\tilde b_n:=b_n+\frac{1}{2}\sigma^{jl}_n\p_j\sigma^{\cdot l}_n$.
We divide the proof into three steps.

{\bf (Step 1)}.  By Lemma \ref{Le5} and the property of the convolution, for all $x\in\mR^d$ and $t\in[0,1]$, we have
\begin{align*}
\mE|\det(\nabla [X^n_t(x)]^{-1})|&\leq\exp\Big\{\|[-\div \tilde b_n+\tfrac{1}{2}\p_i\sigma^{jl}_n\p_j\sigma^{il}_n
+\sigma^{il}_n\p_{ij}^2\sigma^{jl}_n+\tfrac{1}{2}|\div \sigma_n|^2]^+\|_\infty\Big\}\\
&=\exp\Big\{\|[-\div  b_n+\tfrac{1}{2}\sigma^{il}_n\p_{ij}^2\sigma^{jl}_n+\tfrac{1}{2}|\div \sigma_n|^2]^+\|_\infty\Big\}\\
&\leq\exp\Big\{\|[\div b_n]^-\|_\infty+\tfrac{1}{2}\||\sigma_n|\cdot|\nabla\div\sigma_n|\|_\infty
+\tfrac{1}{2}\|\div \sigma_n\|^2_\infty\Big\}.
\end{align*}
Noticing that
$$
\div b_n=\p_i\chi_n (b^i*\rho_n)+(\div b*\rho_n)\chi_n,
$$
$$
\sigma^{il}_n\p_{ij}^2\sigma^{jl}_n=(\sigma^{ij}*\rho_n)[(\p^2_{ij}\sigma*\rho_n)\chi_n
+2(\p_i\sigma*\rho_n)\p_j\chi_n+(\sigma*\rho_n)\p^2_{ij}\chi_n],
$$
by (\ref{Con2}), the definition of $\chi_n$ and elementary calculus, for $n>2(\frac{1}{\eps}\vee r)$,
where $r$ is from (\ref{Con2}), we find
\begin{align*}
\|[\div b_n]^-\|_\infty&\leq C+\|[\div b]^-\|_\infty,\\
\||\sigma_n|\cdot|\nabla\div\sigma_n|\|_\infty&\leq
C+\big\|\sup_{|z|\leq \eps}|\sigma(\cdot-z)|\cdot|\nabla\div\sigma|\big\|_\infty,\\
\|\div \sigma_n\|^2_\infty&\leq C+\|\div \sigma\|^2_\infty.
\end{align*}
Here and below,  $C$ is independent of $n$. Thus,
\begin{align*}
\sup_{n\in\mN}\sup_{(t,x)\in[0,1]\times\mR^d}\mE|\det(\nabla [X^n_t(x)]^{-1})|<+\infty.
\end{align*}
Hence, for any nonnegative measurable function $\varphi\in L^1(\mR^d)$,
\begin{align}
\sup_{t\in[0,1]}\mE\int_{\mR^d}\varphi(X^n_t(x))\dif x
=\sup_{t\in[0,1]}\mE\int_{\mR^d}\varphi(x)\cdot|\det(\nabla [X^n_t(x)]^{-1})|\dif x
\leq K\|\varphi\|_{L^1}.\label{PP1}
\end{align}

{\bf (Step 2)}. In this step we prove that for any $N>0$,
\begin{align}
\sup_{n\in\mN}\mE\int_{B_N}\sup_{t\in[0,1]}|X^n_t(x)|^2\dif x<+\infty.\label{PP5}
\end{align}
Set
$$
g_t(x):=\mE\left(\sup_{s\in[0,t]}|X^n_s(x)|^2\right).
$$
By It\^o's formula, Burkholder's inequality and Young's inequality, we have
\begin{align*}
g_t(x)&\leq |x|^2+2\mE\int^t_0|X^n_s(x)|\cdot|\tilde b_n(X^n_s(x))|\dif s+\mE\int^t_0\|\sigma_n(X^n_s(x))\|^2\dif s\\
&\quad+C\mE\left(\int^t_0|X^n_s(x)|^2\cdot\|\sigma_n(X^n_s(x))\|^2\dif s\right)^{1/2}\\
&\leq |x|^2+2\mE\int^t_0|X^n_s(x)|\cdot|\tilde b_n(X^n_s(x))|\cdot (1_{|X^n_s(x)|\leq r}+1_{|X^n_s(x)|>r})\dif s\\
&\quad+\mE\int^t_0\|\sigma_n(X^n_s(x))\|^2\dif s+C\mE\left(\sup_{s\in[0,t]}|X^n_s(x)|\left[\int^t_0\|\sigma_n(X^n_s(x))\|^2\dif s\right]^{1/2}\right)\\
&\leq |x|^2+2r\mE\int^t_0|\tilde b_n(X^n_s(x))|\cdot 1_{|X^n_s(x)|\leq r}\dif s
+C_r\mE\int^t_0(1+|X^n_s(x)|^2)\dif s\\
&\quad+\frac{1}{2}g_t(x)+C\mE\int^t_0\|\sigma_n(X^n_s(x))\|^2\dif s,
\end{align*}
where $r$ is from (\ref{Con2}) and we have used (\ref{Con2}) in the last step. Hence,
\begin{align*}
g_t(x)&\leq 2|x|^2+4r\mE\int^t_0|\tilde b_n(X^n_s(x))|\cdot 1_{|X^n_s(x)|\leq r}\dif s\\
&\quad +2C_r\int^t_0(1+g_s(x))\dif s+C\mE\int^t_0\|\sigma_n(X^n_s(x))\|^2\dif s.
\end{align*}
By Gronwall's inequality, we obtain that
$$
g_1(x)\leq C_r\left(|x|^2+\mE\int^1_0|\tilde b_n(X^n_s(x))|\cdot 1_{|X^n_s(x)|\leq r}\dif s
+\mE\int^1_0\|\sigma_n(X^n_s(x))\|^2\dif s\right).
$$
Now, by (\ref{PP1}) and (\ref{Con2}), we have
\begin{align*}
\mE\int_{B_N}g_t(x)\dif x&\leq C_{N,r}+C_r\|\tilde b_n\|_{L^1(B_r)}
+C_{N,r}(\|\sigma_n\|_{L^\infty(B_r^c)}^2+\|\sigma_n\|^2_{L^2(B_r)})\\
&\leq C_{N,r}+C_r\|b_n\|_{L^1(B_r)}+C_r\|\sigma_n\|_{L^{2}(B_r)}\|\nabla\sigma_n\|_{L^{2}(B_r)}
+C_{N,r}(\|\sigma\|_{L^\infty(B_r^c)}^2+\|\sigma\|^2_{L^2(B_r)})\\
&\leq C_{N,r}+C_r\|b\|_{L^1(B_r)}+C_r\|\sigma\|_{L^{2}(B_r)}\|\nabla\sigma\|_{L^{2}(B_r)}
++C_{N,r}(\|\sigma\|_{L^\infty(B_r^c)}^2+\|\sigma\|^2_{L^2(B_r)}),
\end{align*}
which gives (\ref{PP5}).

{\bf (Step 3)}.
Noting that for $n>R+1$
$$
\|\nabla b_n\|_{L^1(B_{R+1})}\leq \|\nabla b\|_{L^1(B_{R+1})},\ \ \|\nabla\sigma_n\|_{L^2(B_{R+1})}\leq
\|\nabla\sigma\|_{L^2(B_{R+1})},
$$
by (\ref{PP1}), (\ref{PP5}) and Lemma \ref{Le1}, we have that for any $\delta,\eta,\eps\in(0,1)$,
\begin{align*}
\mE\int_{B_N}\left(\sup_{t\in[0,1]}|X^n_t(x)-X^m_t(x)|^{2}\wedge 1\right)& \dif x\leq \eta
+\frac{C(N,r)}{R\eta}+\frac{C_2}{\eta}
\|\nabla(b_n*\varrho_\eps-b_n)\|_{L^1(B_{R+1})}+\frac{C_1\eps^{-d}}{\eta\log\delta^{-1}}\\
&+\frac{C_3}{\eta\delta\log\delta^{-1}}\left(\|b_n- b_m\|_{L^1(B_{R})}
+\|\sigma_n-\sigma_m\|_{L^{2}(B_{R})}\right),
\end{align*}
where $C_1,C_2, C_3$ are independent of $n,\eps,\delta$.

We take limits according to the following order: $n,m\to\infty$, $\delta\to 0$, $\eps\to 0$, $R\to\infty$,  $\eta\to 0$,
then find
$$
\lim_{n,m\to\infty}\mE\int_{B_N}\left(\sup_{t\in[0,1]}|X^n_t(x)-X^m_t(x)|^{2}\wedge 1\right)\dif x=0,
$$
which together with (\ref{PP5}) gives further that for any $p\in[1,2)$,
$$
\lim_{n,m\to\infty}\mE\int_{B_N}\left(\sup_{t\in[0,1]}|X^n_t(x)-X^m_t(x)|^p\right)\dif x=0.
$$
Therefore, there exists a continuous $\sF_t$-adapted stochastic field  $X_t(x)$ such that for any $N>0$ and $p\in[1,2)$,
$$
\lim_{n\to\infty}\mE\int_{B_N}\left(\sup_{t\in[0,1]}|X^n_t(x)-X_t(x)|^p\right)\dif x=0.
$$
In particular, there exists a subsequence still denoted by $n$ such that for $P\otimes\mu$-almost all $(\omega,x)$,
$$
\lim_{n\to\infty}\sup_{t\in[0,1]}|X^{n}_t(\omega,x)-X_t(\omega,x)|=0.
$$
Condition {\bf (A)} in Definition \ref{Def1} now follows by (\ref{PP1}) and (i) of Lemma \ref{Le4}.
For verifying {\bf (B)} in Definition \ref{Def1}, it suffices to prove that for any $N>0$ and $s\in[0,1]$,
\begin{align}
&\lim_{n\to\infty}\mE\int_{B_N}|b_n(X^n_s(x))-b(X_s(x))|\dif x=0,\label{PP2}\\
&\lim_{n\to\infty}\mE\int_{B_N}|(\sigma^{jl}_n\p_j\sigma^{il}_n)(X^n_s(x))-(\sigma^{jl}\p_j\sigma^{il})(X_s(x))|\dif x=0,\label{PP3}\\
&\lim_{n\to\infty}\mE\int_{B_N}|\sigma_n(X^n_s(x))-\sigma(X_s(x))|^2\dif x=0.\label{PP4}
\end{align}
We only prove (\ref{PP2}). The others are analogous. We make the following decomposition:
\begin{align*}
\int_{B_N}\!\!\!|b_n(X^n_s(x))-b(X_s(x))|\dif x&\leq \int_{B_N}\!\!\!|b_n\chi_m-b\chi_m|(X^n_s(x))|\dif x
+\int_{B_N}\!\!\!|b_n(1-\chi_m)|(X^n_s(x))\dif x\\
&\quad+\int_{B_N}\!\!\!|b(1-\chi_m)|(X_s(x))\dif x=:I^{nm}_1+I^{nm}_2+I^{m}_3.
\end{align*}
For fixed $m\in\mN$, by (ii) of Lemma \ref{Le4}, we have
\begin{align}
\lim_{n\to\infty}\mE I^{nm}_1=0.\label{PI1}
\end{align}
On the other hand,  for $m>r$, we have
$$
I^{nm}_2\leq C\int_{B_N}\!\!\!(1+|X^n_s(x)|)\cdot 1_{|X^n_s(x)|\geq m}\dif x
\leq \frac{C}{m}\int_{B_N}\!\!\!(1+|X^n_s(x)|^2)\dif x,
$$
which together with  (\ref{PP5}) yields
\begin{align}
\lim_{m\to\infty}\sup_{n}\mE I^{nm}_2=0.\label{PI2}
\end{align}
Similarly,
\begin{align}
\lim_{m\to\infty}\mE I^{m}_3=0.\label{PI3}
\end{align}
Combining (\ref{PI1}), (\ref{PI2}) and (\ref{PI3}), we get (\ref{PP2}). The proof is thus complete.
\end{proof}

\begin{proof}[Proof of Theorem \ref{Main1}] Let $b_n$ and $\sigma_n$ be defined by (\ref{Es3}).
Since $b$ and $\sigma$ are linear growth, we have
$$
|b_n(x)|+|\sigma_n(x)|\leq C(1+|x|),
$$
where $C$ is independent of $n$. It is then standard to prove that for any $p\geq 1$,
$$
\sup_{n\in\mN}\mE\left(\sup_{t\in[0,1]}|X^n_t(x)|^{2p}\right)<+\infty.
$$
Note that
$$
\p_j\sigma^{il}_n=\p_j\sigma^{il}*\rho_n\cdot\chi_n+\sigma^{il}_n\cdot \p_j\chi_n,
$$
and by the linear growth of $\sigma$
$$
|\sigma_n\cdot\nabla\chi_n|\leq\frac{C1_{n\leq|x|\leq 2n}}{n}\int_{\mR^d}(1+|x-y|)\rho_n(y)\dif y\leq C.
$$
By Jensen's inequality and (\ref{Ep0}), for $n\geq\frac{1}{\eps}$, we have
\begin{align*}
|\Lambda^{\sigma_n}_1|^2&=|\div\sigma_n+\sigma^{i\cdot}_n\p_i\lambda|^2\\
&\leq C\left(|\div\sigma|^2*\rho_n
+|\sigma|^2*\rho_n\cdot|\nabla\lambda|^2+1\right)\\
&\leq C\left(|\nabla\sigma|^2+|\sigma|^2\gamma_2^2\right)*\rho_n+C
\end{align*}
and
\begin{align*}
-\Lambda^{b_n,\sigma_n}_2&=-\Big[\div b_n+b^i_n\p_i\lambda+
\frac{1}{2}(\sigma^{il}_n\sigma^{jl}_n\p^2_{ij}
\lambda-\p_i\sigma^{jl}_n\p_j\sigma^{il}_n)\Big]\\
&\leq C\Big[[\div b]^-*\rho_n+|b|*\rho_n\cdot|\nabla\lambda|+
(|\sigma|*\rho_n)^2\cdot|\nabla^2\lambda|+(|\nabla\sigma|*\rho_n)^2)+1\Big]\\
&\leq C\Big[[\div b]^-+|b|\gamma_2+|\sigma|^2\gamma_3+|\nabla\sigma|^2\Big]*\rho_n+C.
\end{align*}
Hence, for all $t\in[0,1]$ and $p>1$, by Lemma \ref{Le00} and Jensen's inequality again,
\begin{align*}
\mE\int_{\mR^d}|\cJ^n_t(x)|^p\mu(\dif x)&\leq
C_N\sup_{t\in[0,1]}\int_{\mR^d}\exp\Big\{tp^3|\Lambda^{\sigma_n}_3(x)|^2
-tp^2 \Lambda^{b_n,\sigma_n}_2(x)\Big\}\mu(\dif x)\\
& \leq C_N\int_{\mR^d}e^{ C\big([\div b]^-+|b|\gamma_2+|\sigma|^2(\gamma^2_2+\gamma_3)
+|\nabla\sigma|^2\big)*\rho_n(x) }\cdot e^{\lambda(x)}\dif x\\
& \leq C_N\int_{\mR^d}e^{ \big[C\big([\div b]^-+|b|\gamma_2+|\sigma|^2(\gamma^2_2+\gamma_3)
+|\nabla\sigma|^2\big)+\gamma_1\big]*\rho_n(x) }\dif x\\
& \leq C_N\int_{\mR^d}e^{ C\big([\div b]^-+|b|\gamma_2+|\sigma|^2(\gamma^2_2+\gamma_3)
+|\nabla\sigma|^2\big)+\gamma_1} *\rho_n(x)\dif x\\
& = C_N\int_{\mR^d}e^{\big[ C\big([\div b]^-+|b|\gamma_2+|\sigma|^2(\gamma^2_2+\gamma_3)
+|\nabla\sigma|^2\big)+\gamma_1\big](x) }\dif x<+\infty.
\end{align*}
Thus, by (\ref{P0}) and H\"older's inequality, we obtain that for any $p>1$,
\begin{align*}
\mE\int_{\mR^d}&\varphi(X^n_t(x))\mu(\dif x)=
\mE\int_{\mR^d}\varphi(x) \cJ^n_t(x)\mu(\dif x) \leq\|\varphi\|_{L^p_\mu}\left(\mE\int_{\mR^d}
|\cJ^n_t(x)|^{\frac{p}{p-1}}\mu(\dif x)\right)^{1-\frac{1}{p}}\leq C.
\end{align*}
The rest proof is the same as the {\bf Step 3} in the proof of Theorem \ref{Main2}.
\end{proof}

\section{Proof of Theorem \ref{Main3}}

For proving Theorem \ref{Main3}, our task is to check {\bf (LD)$_1$} and {\bf (LD)$_2$}.
By the infinite-dimensional Yamada-Watanabe theorem (cf. \cite{Ro-Sc-Zh}), there exists a measurable functional
$$
\Phi_\eps: \Omega\to\mS=L^{2p}_\nu(\mR^d; C([0,1];\mR^d)),\ \ p\geq 1,
$$
such that
$$
X_{\eps,t}(\omega,x)=\Phi_\eps(\omega)(t,x).
$$
For $\eps\in(0,1)$, let $h^\eps\in\cA_M$, where $\cA_M$ is defined by (\ref{Op2}). By Girsanov's theorem, one sees that
$$
X^\eps_t(\omega,x)=\Phi_\eps\left(W_\cdot(\omega)
+\frac{1}{\sqrt{\eps}}\int^\cdot_0h^\eps_s(\omega)\dif s\right)(t,x)
$$
solves the controlled equation:
$$
\dif X^\eps_t(x)=b(X^\eps_t(x))\dif t+\sigma(X^\eps_t(x))h^\eps_t\dif t+\sqrt{\eps}\sigma(X^\eps_t(x))\dif W_t,
\ \ X^\eps_0(x)=x.
$$
For $h\in\cA_M$, let $X_t^h(x)$ solve  equation (\ref{C1}).
We have:
\bl\label{Le6}
(i). For any $p\geq 1$ and $h\in\cA_M$,
$$
\mE\left(\sup_{t\in[0,1]}|X^h_t(x)|^{2p}\right)+\sup_{\eps\in(0,1)}\mE\left(\sup_{t\in[0,1]}|X^\eps_t(x)|^{2p}\right)\leq C(1+|x|^{2p}).
$$
(ii). For any $p>1$, $h^\eps\in\cA_M$ and nonnegative function $\varphi\in L^p_\mu(\mR^d)$,
$$
\mE\int_{B_N}\varphi(X^\eps_t(x))\mu(\dif x)\leq C_{N,M}\|\varphi\|_{L^p_\mu}.
$$
\el
\begin{proof}
(i). It is standard by the linear growth of $b$ and $\sigma$.

(ii). Let us define $b_n$ and $\sigma_n$ by (\ref{Es3}). Consider the following SDE:
$$
\dif X^{\eps,n}_t(x)=b_n(X^{\eps,n}_t(x))\dif t
+\sigma_n(X^{\eps,n}_t(x))h^\eps_t\dif t+\sqrt{\eps}\sigma_n(X^{\eps,n}_t(x))\dif W_t,\ \ X^{\eps,n}_0(x)=x.
$$
From the proof of Lemma \ref{Le00} and Theorem \ref{Main1}, one can see that for any $p>1$ and $\varphi\in L^p_\mu(\mR^d)$,
\begin{align*}
\mE\int_{B_N}\varphi(X^{\eps,n}_t(x))\mu(\dif x)
\leq C_{N,M}\|\varphi\|_{L^{p}_\mu},
\end{align*}
where $C_{N,M}$ is independent of $\eps$. Now taking limit $n\to\infty$ gives the result (see Lemma \ref{Le4}).
\end{proof}

Set
\begin{align}
w^\eps_t(x):=\int^t_0\sigma(X^h_s(x))(h^\eps_s-h_s)\dif s.\label{PP6}
\end{align}
\bl\label{Le9}
Suppose that $h_\eps$ weakly converges to $h$ a.s.  in $\cD_M$. Then for any $p\geq 1$, we have
$$
\lim_{\eps\to 0}\mE\int_{B_N}\sup_{t\in[0,1]}|w^\eps_t(x)|^{2p}\dif x=0.
$$
\el
\begin{proof}
For fixed $(\omega,x)$, let us first prove that
\begin{align}
\lim_{\eps\to 0}\sup_{t\in[0,1]}|w^\eps_t(\omega,x)|=0.\label{Es4}
\end{align}
By the weak convergence of $h^\eps_\cdot(\omega)$ to $h_\cdot(\omega)$, one sees that for fixed $t\in[0,1]$
$$
\lim_{\eps\to 0}w^\eps_t(\omega,x)=\lim_{\eps\to 0} \int^t_0\sigma(X^h_s(\omega,x))(h^\eps_s(\omega)-h_s(\omega))\dif s=0.
$$
Since for $t'<t$
\begin{align*}
|w^\eps_t(\omega,x)-w^\eps_{t'}(\omega,x)|&\leq\int^t_{t'}|\sigma(X^h_s(\omega,x))(h^\eps_s(\omega)-h_s(\omega))|\dif s\\
&\leq2M\left(\int^t_{t'}|\sigma(X^h_s(\omega,x))|^2\dif s\right)^{\frac{1}{2}}\to 0,
\end{align*}
uniformly in $\eps$ as $|t-t'|\to 0$,
we immediately have (\ref{Es4}). In view of
$$
\sup_{t\in[0,1]}|w^\eps_t(x)|^{2p}\leq C_{M,p}\int^1_0|\sigma(X^h_s(x))|^{2p}\dif s,
$$
the desired limit now follows by the dominated convergence theorem and (\ref{Es4}).
\end{proof}

\bl\label{Le10}
Suppose that $h^\eps$ weakly converges to $h$ a.s.  in $\cD_M$.
Then for some subsequence $\eps_k$, $X^{\eps_k}$ converges to $X^h$ in probability in space $\mS$,
where $X^h$ solves  equation (\ref{C1}).
\el
\begin{proof}
Set
$$
Z^\eps_t(x):=X^\eps_t(x)-X^h_t(x).
$$
By It\^o's formula, for any $\delta>0$, we have
\begin{align*}
\log\left(\frac{|Z^\eps_t(x)|^2}{\delta^2}+1\right)
&=2\int^t_0\frac{\<Z^\eps_s(x), b(X^\eps_s(x))-b(X^h_s(x))\>}{|Z^\eps_s(x)|^2+\delta^2}\dif s\\
&\quad+2\int^t_0\frac{\<Z^\eps_s(x), (\sigma(X^\eps_s(x))-\sigma(X^h_s(x)))h^\eps_s\>}{|Z^\eps_s(x)|^2+\delta^2}\dif s\\
&\quad+2\int^t_0\frac{\<Z^\eps_s(x), \sigma(X^h_s(x))(h^\eps_s-h_s)\>}{|Z^\eps_s(x)|^2+\delta^2}\dif s\\
&\quad+2\sqrt{\eps}\int^t_0\frac{\<Z^\eps_s(x), \sigma(X^\eps_s(x))\dif W_s\>}{|Z^\eps_s(x)|^2+\delta^2}\\
&\quad+\eps\int^t_0\frac{\|\sigma(X^\eps_s(x))\|^2}{|Z^\eps_s(x)|^2+\delta^2}\dif s
-2\eps\int^t_0\frac{|(\sigma(X^\eps_s(x)))^\mathrm{t}\cdot Z^\eps_s(x)|^2}{(|Z^\eps_s(x)|^2+\delta^2)^2}\dif s\\
&=:I^\eps_1(t,x)+I^\eps_2(t,x)+I^\eps_3(t,x)+I^\eps_4(t,x)+I^\eps_5(t,x)+I^\eps_6(t,x).
\end{align*}
We want to prove that for any $N, R>0$,
\begin{align}
\mE\int_{B_N\cap G^\eps_R}\log\left(\frac{\sup_{t\in[0,1]}|Z^\eps_t(x)|^2}{\delta^2}+1\right)\mu(\dif x)\leq
C_1+\frac{C_2(\eps)}{\delta},\label{Po2}
\end{align}
where $C_1$ is independent of $\eps$ and $\delta$, $C_2(\eps)\to 0$ as $\eps\to 0$ and
$$
G^\eps_R(\omega):=\Big\{x\in\mR^d: \sup_{t\in[0,1]}|X^\eps_t(\omega,x)|\vee|X^h_t(\omega,x)|\leq R\Big\}.
$$
First of all, $I^\eps_6(t,x)$ is negative and dropped. By Lemmas \ref{Le0} and \ref{Le6},
as in the proof of Lemma \ref{Le1}, it is easy to see that
$$
\mE\int_{B_N\cap G^\eps_R}\sup_{t\in[0,1]}(|I^\eps_1(t,x)|+I^\eps_2(t,x)|)\mu(\dif x)\leq C_1.
$$
Moreover, by Burkholder's inequality, we also have
$$
\mE\int_{B_N\cap G^\eps_R}\sup_{t\in[0,1]}(|I^\eps_4(t,x)|+I^\eps_5(t,x)|)\mu(\dif x)\leq \frac{C\eps}{\delta^2}.
$$

We now deal with the hard term $I^\eps_3(t,x)$. Set
$$
\xi(x):=\frac{x}{|x|^2+\delta^2}.
$$
Recalling (\ref{PP6}), we have
$$
I^\eps_3(t,x)=2\int^t_0\<\xi(Z^\eps_s(x)), \dif w^\eps_s(x)\>
=2\<\xi(Z^\eps_t(x)), w^\eps_t(x)\>-2\int^t_0\<w^\eps_s(x),\dif \xi(Z^\eps_s(x))\>.
$$
By It\^o's formula, we have
\begin{align*}
\dif \xi(Z^\eps_t(x))&=\nabla \xi(Z^\eps_t(x))(b(X^\eps_t(x))-b(X^h_t(x)))\dif t\\
&\quad+\nabla \xi(Z^\eps_t(x))(\sigma(X^\eps_t(x))h^\eps_t-\sigma(X^h_t(x))h_t)\dif t\\
&\quad+\frac{\eps}{2}\p^2_{ij}\xi(Z^\eps_t(x))\sigma^{il}(X^\eps_t(x))\sigma^{jl}(X^\eps_t(x))\dif t\\
&\quad+\sqrt{\eps}\nabla \xi(Z^\eps_t(x))\sigma(X^\eps_t(x))\dif W_t.
\end{align*}
Hence,
\begin{align*}
I^\eps_3(t,x)&=2\<\xi(Z^\eps_t(x)), w^\eps_t(x)\>-2\int^t_0\<\nabla \xi(Z^\eps_s(x))
(b(X^\eps_s(x))-b(X^h_s(x))), w^\eps_s(x)\>\dif s\\
&\quad-2\int^t_0\<\nabla \xi(Z^\eps_s(x))(\sigma(X^\eps_s(x))h^\eps_s-\sigma(X^h_s(x))h_s), w^\eps_s(x)\>\dif s\\
&\quad-\eps\int^t_0\<\p^2_{ij}\xi(Z^\eps_s(x))\sigma^{il}(X^\eps_s(x))\sigma^{jl}(X^\eps_s(x)), w^\eps_s(x)\>\dif s\\
&\quad-2\sqrt{\eps}\int^t_0\<\nabla \xi(Z^\eps_s(x))\sigma(X^\eps_s(x))\dif W_s, w^\eps_s(x)\>\\
&=:I^\eps_{31}(t,x)+I^\eps_{32}(t,x)+I^\eps_{33}(t,x)+I^\eps_{34}(t,x)+I^\eps_{35}(t,x).
\end{align*}
Noticing that
$$
\p_i \xi^k(x)=\frac{1_{i=k}}{|x|^2+\delta^2}-\frac{2x^ix^k}{(|x|^2+\delta^2)^2}
$$
and
$$
\p^2_{ij} \xi^k(x)=-\frac{2\cdot 1_{i=k}x^j}{(|x|^2+\delta^2)^2}+\frac{4x^ix^jx^k}{(|x|^2+\delta^2)^3},
$$
we have
$$
|\xi(x)|\leq\frac{1}{\delta},\ \ |\nabla \xi(x)|\leq\frac{2}{\delta^2},\ \ |\nabla^2 \xi(x)|\leq\frac{6}{\delta^3}.
$$
Using Lemma \ref{Le9}, as above, one finds that
$$
\mE\int_{B_N\cap G^\eps_R}\sup_{t\in[0,1]}|I^\eps_3(t,x)| \mu(\dif x)\leq \frac{C(\eps)}{\delta^3},
$$
where $C(\eps)\to 0$ as $\eps\to 0$.

Combining the above estimates, we obtain (\ref{Po2}). Thus,
by (\ref{Po2}) and Lemma \ref{Le6}, as {\bf (Step 2)} in  the proof of Lemma \ref{Le1},
there exists a subsequence $\eps_k$ such that for $P\otimes\mu$-almost all $(\omega,x)$
$$
\sup_{t\in[0,1]}|X^{\eps_k}_t(\omega,x)-X^h_t(\omega,x)|\to 0,\ \ \mbox{ as $k\to\infty$.}
$$
Using (i) of Lemma \ref{Le6}, there exists another subsequence $\eps_k'$ such that $X^{\eps'_k}$
converges to $X^h$ in probability in space $\mS$.
\end{proof}

\begin{proof}[Proof of Theorem \ref{Main3}]
Let $h^\eps$ be a sequence in $\cA_M$ converging to $h$ in distribution.
Since $\cD_M$ is compact and the law of $W$ is tight,
$\{h^\eps,W\}$ is tight in $\cD_M\times\Omega$ by the definition of tightness.
Without loss of generality, we assume that the law of $\{h^\eps,W\}$ weakly converges to some $\mP$ on $\cD_M\times\Omega$.
Then the law of $h$ is just $\mP(\cdot,\Omega)$.
By Skorokhod's representation theorem, there are probability space $(\tilde \Omega,\tilde\sF, \tilde P)$,  and random elelments
$\{\tilde h^\eps,\tilde W^\eps\}$ and $\{\tilde h,\tilde W\}$ in $\cD_M\times\Omega$ such that

(1) $(\tilde h^\eps,\tilde W^\eps)$ a.s. converges to $(\tilde h,\tilde W)$;

(2) $(\tilde h^\eps,\tilde W^\eps)$ has the same law as $(h^\eps,W)$;

(3) The law of $\{\tilde h,\tilde W\}$ is $\mP$, and the law of
$h$ is the same as $\tilde h$.

Using Lemma \ref{Le10}, we get for some subsequence $\eps_k$,
$$
\Phi_{\eps_k}\left(\tilde W^{\eps_k}_\cdot+\frac{1}{\sqrt{\eps_k}}\int^\cdot_0
\tilde h^{\eps_k}_s\dif s\right)\to X^{\tilde h}, \ \ \mbox{ in
probability}.
$$
From this, we derive
$$
\Phi_{\eps_k}\left(W_\cdot+\frac{1}{\sqrt{\eps_k}}\int^\cdot_0h^{\eps_k}_s\dif s\right)\to X^{h}, \ \mbox{ in distribution}.
$$
Thus, {\bf (LD)$_\mathbf{1}$} holds.  {\bf (LD)$_\mathbf{2}$} can be simply verified as Lemma \ref{Le10}.
\end{proof}

{\bf Acknowledgements:}

The author is very grateful to Professor Nicolas Privault for providing him an opportunity to
work in the City University of Hong Kong. This work was done during his very pleasant stay in Hong Kong.
The supports by NSFs of China (Nos. 10971076; 10871215) are also acknowledged.


\begin{thebibliography}{999}

\bibitem{Am}Ambrosio, L.: Transport equation and
Cauchy problem for $BV$ vector fields.  Invent. Math.,  158  (2004),  no. 2, 227--260.

\bibitem{Am-Fi}Ambrosio, L.; Figalli, A.: On flows associated to Sobolev vector fields in Wiener spaces:
an approach \'a la DiPerna-Lions.  J. Funct. Anal.  256  (2009),  no. 1, 179--214.

\bibitem{Bo-Du} Bou\'e, M. and Dupuis, P.: A variational representation for certain
functionals of Brownian motion. Ann. of Prob., Vol. 26, No. 4, 1641-1659 (1998).

\bibitem{Bo-Du-El}Bou\'e, M.; Dupuis, P.; Ellis, R. S.: Large deviations for small noise
diffusions with discontinuous statistics.  Probab. Theory Related Fields  116  (2000),  no. 1, 125--149.

\bibitem{bd0}Budhiraja, A. and Dupuis, P.: A variational representation
for positive functionals of infinite dimensional Brownian motion.
Probab. Math. Statist., 20 (2000), no. 1,
Acta Univ. Wratislav. No. 2246, 39--61.

\bibitem{Bu-Du-Ma}Budhiraja, A., Dupuis,P.  and Maroulas, V.: Large deviations for infinite dimensional
stochastic dynamical systems.  Ann. of Prob.,  36 ,  no. 4, 1390--1420 (2008).

\bibitem{Cr}Cruzeiro, A.B.: \'Equations diff¨¦rentielles ordinaires: non explosion et mesures quasi-invariantes.
J. Funct. Anal.  54  (1983),  no. 2, 193--205.

\bibitem{Ci-Cr}Cipriano, F., Cruzeiro, A.B.: Flows associated with irregular $\mR^d$-vector fields.
 J. Differential Equations  219  (2005),  no. 1, 183--201.

\bibitem{Cr-De-Le}Crippa G. and De Lellis C.: Estimates and regularity results for
the DiPerna-Lions flow. J. reine angew. Math. 616 (2008), 15-46.

\bibitem{Di-Li}DiPerna R.J. and Lions P.L.: Ordinary differential equations, transport theory
and Sobolev spaces.  Invent. Math., 98,511-547(1989).

\bibitem{Du-El}Dupuis, P. and Ellis, R.S.: A Weak Convergence Approach to the
Theory of Large Deviations. Wiley, New-York, 1997.

\bibitem{Fa-Lu}Fang, S. and Luo, D.: Transport equations and quasi-invariant flows on the Wiener space.
Bulletin des Sciences Math\'ematiques, doi:10.1016/j.bulsci.2009.01.001.

\bibitem{Fa-Lu-Th}Fang, S., Luo, D. and Thalmaier, A.: Stochastic differential equations with coefficients
in Sobolev spaces. J. Func. Anal., 259:1129-1168.

\bibitem{Fi}Figalli, A.: Existence and uniqueness of martingale solutions for SDEs with
rough or degenerate coefficients.  J. Funct. Anal.  254  (2008),  no. 1, 109--153.

\bibitem{Gy-Ma}Gy\"ongy I., Martinez, T.: On stochastic differential
equations with locally unbounded drift.  Czechoslovak Math. J.
51(126) (2001), no. 4, 763--783.

\bibitem{Ik-Wa}Ikeda, N. and Watanabe S.: Stochastic differential equations
and diffusion processes. North-Holland/Kodanska, Amsterdam/Tokyo, 1981.



\bibitem{Ku}Kunita, H.: Stochastic flows and stochastic differential equations.
Cambridge, Cambridge University Press, 1990.

\bibitem{Kr-Ro}Krylov, N.V. and R\"ockner, M.: Strong solutions of stochasitc equations
with singluar time dependent drift. Probab. Theory Related Fields  131  (2005),  no. 2, 154--196.

\bibitem{Le-Li}Le Bris, C.  and Lions, P.L. : Renormalized solutions of some transport equations
with partially $W^{1,1}$ velocities and applications.
Annali di Matematica, 183, 97-130(2004).

\bibitem{Le-Li1}Le Bris, C.  and Lions, P.L. : Existence and uniqueness of solutions to
Fokker-Planck type equations with irregular coefficients.
Comm. in Partial Diff. Equ., 33:1272-1317,2008.

\bibitem{Ma}Malliavin, P.: Stochastic analysis. Springer-Verlag, 1997.

\bibitem{Re-Xu-Zh}Ren, J., Xu, S. and Zhang, X.: Large deviations for multi-valued stochastic differential equations,
to appear in Journal of Theoretical Probability.

\bibitem{RenJie-Zh}Ren, J. and Zhang, X.: Limit theorems for stochastic differential equations
with discontinuous coefficients. Preprint.

\bibitem{Ro-Zh}R\"ockner, M. and Zhang, X.: Weak uniqueness of Fokker-Planck equations
with degenerate and bounded coefficients.  Comptes Rendus Mathematique, 348, (2010) 435-438.

\bibitem{Ro-Sc-Zh}R\"ockner, M., Schmuland, B., Zhang, X.: Yamada-Watanabe Theorem for
Stochastic Evolution Equations in Infinite Dimensions.  Condensed Matter Physics,
Vol.11, No.2(54), 247-259(2008).

\bibitem{St}Stein, E. M.: Singular integrals and differentiability properties of functions.
Princeton, N.J.,  Princeton University Press,  1970.

\bibitem{Zh3}Zhang, X.: Strong solutions of SDEs with singular drift and
Sobolev diffusion coefficients. Stoch. Proc. and Appl., 115/11 pp. 1805-1818(2005).

\bibitem{Zh2}Zhang, X.: Stochastic flows of SDEs with irregular coefficients and
stochastic transport equations. Bull. Sci. Math. France,  Vol. 134, (2010) 340-378.

\bibitem{Zh5}Zhang, X.: A variational representation for random
functionals on abstract Wiener spaces. J. Math. Kyoto Univ., Vol.49(3), 475-490(2009).

\bibitem{Zv}Zvonkin, A.K.: A transformation of the phase space of a diffusion process
that removes the drift. Mat. Sbornik, No.1, 93(135)(1974).
\end{thebibliography}
\end{document}